\documentclass{amsart}
\usepackage{amssymb}

\newcommand{\Coeff}{\operatorname{Coeff}}

\newcommand{\GL}{\operatorname{GL}}
\newcommand{\gr}{\operatorname{gr}}

\newcommand{\Ker}{\operatorname{Ker}}
\renewcommand{\Im}{\operatorname{Im}}
\newcommand{\Lin}{\operatorname{Lin}}
\newcommand{\SL}{\operatorname{SL}}
\newcommand{\T}{\operatorname{T}}
\newcommand{\tr}{\operatorname{tr}}
\newcommand{\Stab}{\operatorname{Stab}}

\newcommand{\D}{\Delta}

\newcommand{\al}{\alpha}
\newcommand{\be}{\beta}
\newcommand{\ga}{\gamma}
\newcommand{\de}{\delta}
\newcommand{\eps}{\varepsilon}

\newcommand{\ka}{\kappa}
\newcommand{\La}{\Lambda}
\newcommand{\la}{\lambda}
\newcommand{\sig}{\sigma}
\renewcommand{\phi}{\varphi}
\newcommand{\om}{\omega}

\newcommand{\C}{\mathbf C}
\newcommand{\N}{\mathbf N}
\newcommand{\Z}{\mathbf Z}

\renewcommand{\P}{\mathbb P}

\newcommand{\A}{{\mathcal A}}
\newcommand{\B}{{\mathcal B}}
\newcommand{\G}{{\mathcal G}}
\renewcommand{\L}{{\mathcal L}}
\newcommand{\M}{{\mathcal M}}
\newcommand{\MM}{{\mathcal N}}
\renewcommand{\O}{{\mathcal O}}
\newcommand{\Con}{\mathcal C}

\newcommand{\tp}{\,{}^\tau\!}
\renewcommand{\o}{\otimes}
\newcommand{\ol}{\overline}
\newcommand{\inv}{^{-1}}

\newcommand{\pmat}[1]{\begin{pmatrix}#1\end{pmatrix}}
\newcommand{\Pp}{P^{+}}
\newcommand{\rc}{\backslash}
\newcommand{\exch}{\leftrightarrow}
\newcommand{\bR}{R^*}
\newcommand{\bS}{S^*}
\newcommand{\bSp}{\tilde S^*}
\newcommand{\Sym}{\mathrm{Sym}}

\newcommand{\xto}{\xrightarrow}

\newcommand{\spmat}[1]{\left(\begin{smallmatrix}#1\end{smallmatrix}\right)}

\newcommand{\sll}{\operatorname{\mathfrak{sl}}}

\theoremstyle{definition}
\newtheorem{df}{Definition}[section]

\theoremstyle{plain}
\newtheorem{theo}[df]{Theorem}
\newtheorem{prop}[df]{Proposition}
\newtheorem{lem}[df]{Lemma}
\newtheorem{cor}[df]{Corollary}

\numberwithin{equation}{section}

\title[Quantum $\SL(3)$'s]{Quantum $\SL(3,\C)$'s with classical representation theory}
\author{Christian Ohn}
\address{Universit\'e de Reims\\ D\'epartement de Math\'ematiques
(URA CNRS 1870)\\ Moulin de la Housse, B.P. 1039\\ F-51687 Reims cedex 2\\
France}
\email{christian.ohn@univ-reims.fr}

\begin{document}

\begin{abstract}
We study and classify almost all quantum $\SL(3,\C)$'s whose representation
theory is ``similar'' to that of the (ordinary) group $\SL(3,\C)$. Only one
case, related to smooth elliptic curves, could not be treated completely.
\end{abstract}

\maketitle

\section{Introduction}

There are several approaches to the theory of quantum groups, depending on
what aspect of group theory one wants to ``quantize'': one may consider the
fundamental object to be, among other things, the algebra of continuous
functions on a compact group \cite{Wo87}, or the enveloping algebra of a
complex semi-simple Lie algebra \cite{Dr,Ji}, or the algebra of polynomial
functions on a complex algebraic group \cite{FRT}. Let us therefore begin
by saying (definition~\ref{Gqdef} below) what we mean by a quantum analogue
of (the algebra of polynomial functions on) a connected complex reductive
group; roughly speaking, we ask the representation theory to be preserved.
We do not claim that this is the only reasonable way to define ``quantum
reductive groups'' (it is easy to see that it excludes, e.g., quantum
tori), but it will be the one taken throughout the present paper.

So let $G$ be a connected complex reductive group, $B$ a Borel subgroup of
$G$ and $P$ (resp.~$\Pp$) the set of integral (resp.\ dominant integral)
weights of $G$ w.r.t.~$B$. For each $\la\in\Pp$, denote by $L_\la$ the
simple $G$-module of highest weight $\la$ and let $d_\la=\dim L_\la$. If
$\la,\mu,\nu\in\Pp$, let $m_{\la\mu\nu}$ be the multiplicity of $L_\nu$ in
the decomposition of $L_\la\o L_\mu$.
\begin{df}
\label{Gqdef}
We call a \emph{quantum $G$} any (not necessarily commutative) Hopf algebra
$\A$ (over $\C$) such that
\begin{enumerate}
\item there is a family $\{V_\la\mid\la\in\Pp\}$ of simple and pairwise
nonisomorphic $\A$-comodules, with $\dim V_\la=d_\la$,
\item every $\A$-comodule is isomorphic to a direct sum of these,
\item for every $\la,\mu\in\Pp$, $V_\la\o V_\mu$ is isomorphic to
$\bigoplus_{\nu} m_{\la\mu\nu} V_\nu$.
\end{enumerate}
\end{df}
Of course, the algebra $\O(G)$ of polynomial functions on the (ordinary)
group $G$ is a quantum $G$.

Recall that the tensor product of two $\A$-comodules is defined using the
algebra structure on $\A$; thus, the idea is that the way a $V_\nu$ sits
inside $V_\la\o V_\mu$ influences this algebra structure. In particular,
the decomposition of $V_\la\o V_\la$ for a given quantum $G$ may well be
such that the exterior square ${\wedge}^2V_\la$ is not a subcomodule; this
basically accounts for the non-commutativity of such a quantum $G$.

There are (at least) two natural notions of equivalence for quantum $G$'s,
namely
\begin{itemize}
\item isomorphism of Hopf algebras,
\item $\C$-linear monoidal equivalence of categories of comodules (called
\emph{categorial} equivalence).
\end{itemize}
The second notion is weaker than the first. For example, the so-called
``Jordanian'' quantum $\SL(2)$ (introduced in~\cite{DMMZ}) is categorially
equivalent to $\O(\SL(2))$~\cite{Wo91}, but it is not commutative.

Up to categorial equivalence, quantum $\SL(n)$'s have been classified
in~\cite{KW}: they are parametrized by a ``deformation'' parameter
(either~$1$ or not a root of unity) and a ``twisting'' parameter (an $n$-th
root of unity).

Up to Hopf algebra isomorphism, quantum $\SL(2)$'s have been classified
in~\cite{Wo91}.

~

In the present work, we study quantum $\SL(3)$'s, up to Hopf algebra
isomorphism. To read definition~\ref{Gqdef} in this case, recall that
$P=\Z^2$, $\Pp=\N^2$, and that $d_{(k,\ell)}=(k+1)(\ell+1)(k+\ell+2)/2$ for
$(k,\ell)\in\N^2$. Also, the multiplicities
$m_{(k,\ell)\,(k',\ell')\,(k'',\ell'')}$ can be computed combinatorially
(using, e.g., the Littlewood-Richardson rule).

If $\A$ is a quantum $\SL(3)$, the idea is to find data consisting of a
finite number of $\A$-comodules and a finite number of $\A$-comodule
morphisms between tensor products of them, such that $\A$ can be
reconstructed (in the Tannaka-Krein sense) from these data, and to see that
classifying quantum $\SL(3)$'s up to isomorphism amounts to classifying
these finite-dimensional data up to (a suitable notion of) equivalence.

In principle, these data could involve only the ``natural'' $3$-dimensional
comodule $V_{(1,0)}$, because definition~\ref{Gqdef}(c) implies that every
$V_{(k,\ell)}$ is contained in some tensor power of $V_{(1,0)}$, so its
matrix coefficients generate $\A$ as an algebra. However, we rather use
\emph{both} ``fundamental'' comodules $V:=V_{(1,0)}$ and $W:=V_{(0,1)}$
(together with a suitable collection of morphisms): the point is that with
these $18$ generators (instead of $9$), $\A$ can be presented by
\emph{quadratic} relations (instead of cubic ones, such as a ``quantum
determinant'').

In section~\ref{AtoLsection}, we make these finite-dimensional data
precise: starting from a given quantum $\SL(3)$, we choose eight morphisms
between tensor products of its comodules $V$ and $W$. The Schur lemma
imposes some compatibility conditions between these morphisms. This leads
us to the definition of a \emph{basic quantum datum} (BQD for short), in
which $V$ and $W$ become just vector spaces and the eight morphisms just
linear maps, satisfying the compatibility conditions mentioned above.

Conversely, in section~\ref{LtoAsection}, we start from a BQD $\L$ and we
reconstruct a Hopf algebra $\A_\L$ by the usual Tannaka-Krein procedure.
The goal of sections~\ref{shapesection}--\ref{simplesection} is to see
whether this Hopf algebra is actually a quantum $\SL(3)$. In other words:
if the fundamental comodules of a Hopf algebra are ``$\SL(3)$-ish'', does
it follow that \emph{all} comodules are?

To understand sections~\ref{shapesection} and~\ref{filtersection}, let us
first recall the following well-known situation. Let $G,B,P,\Pp$ be as at
the beginning of this introduction, denote by $U$ the unipotent radical of
$B$ and let $T$ be a maximal torus in $B$. View elements of $P$ as
characters of $T$. Since $T$ normalizes $U$, $T$ acts from the right on
$G/U$ and from the left on $U\rc G$; this induces $\Pp$-gradings
$\O(\ol{G/U})=\bigoplus_{\la\in\Pp}V_\la$ and $\O(\ol{U\rc
G})=\bigoplus_{\la\in\Pp}V^\la$. By the Borel-Weil theorem, the $V_\la$'s
are precisely the irreducible representations of $G$; therefore
$\O(\ol{G/U})$ is called a \emph{shape algebra} for $G$. Furthermore,
$V^\la$ can be identified with the dual of $V_\la$, so the algebra of
$T$-invariants
\[
\G(G):=\O(\ol{U\rc G}\times \ol{G/U})^T=\bigoplus_{\la\in\Pp}V^\la\o V_\la
\]
identifies with $\O(G)$ as a vector space, by the Peter-Weyl decomposition.
Actually, for a suitable additive function $h:\Pp\to\N$, $\O(G)$ becomes
$\N$-filtered by putting $V^\la\o V_\la$ into degree~$h(\la)$, and then
$\G(G)\simeq\gr\O(G)$.

(Note that we have avoided using the opposite unipotent subgroup $U^{-}$.
Also, the maximal torus $T$ only appears in the guise of a $P$-grading and
is not really used as a subgroup of $G$. This is necessary, because there
exist quantum $G$'s in which neither $G/U^{-}$ nor $T$ have quantum
analogues: the Jordanian quantum $\SL(2)$ already mentioned is an easy
example.)

In section~\ref{shapesection}, we define two $\N^2$-graded quadratic
algebras $\M_\L=\bigoplus V_{(k,\ell)}$ and $\MM_\L=\bigoplus V^{(k,\ell)}$
(generated by $V\oplus W$ and $V^*\oplus W^*$, respectively), which are
quantum analogues of $\O(\ol{G/U})$ and $\O(\ol{U\rc G})$ (for $G=\SL(3)$).
We show that $\dim V_{(k,\ell)}=\dim V^{(k,\ell)}=d_{(k,\ell)}$ for all
$(k,\ell)$ and that $\M_\L$, $\MM_\L$ are Koszul algebras, except possibly
when $\L$ is a so-called \emph{elliptic} BQD (case I.h in the
classification of section~\ref{classsection}).

In section~\ref{filtersection}, we consider the subalgebra $\G_\L$ of
$\MM_\L\o\M_\L$ defined by
\[
\G_\L:=\bigoplus_{(k,\ell)\in\N^2}V^{(k,\ell)}\o V_{(k,\ell)}
\]
(a quantum analogue of $\G(G)$), which is also $\N^2$-graded by putting
$V^{(k,\ell)}\o V_{(k,\ell)}$ into degree $(k,\ell)$. We give a
presentation of $\G_\L$ that can be deduced from a suitable presentation of
$\A_\L$ by ``cutting off'' all terms of degree $<2$; this yields a
canonical surjection $\G_\L\to\gr\A_\L$. If $\L$ is not elliptic, we show
(noting that $\G_\L$ is still Koszul) that this surjection is an
isomorphism (a quantum analogue of $\G(G)\simeq\gr\O(G)$), using the
results of~\cite{BG}. Here lies the main advantage in dealing with algebras
that are quadratic.

Section~\ref{Heckesection} is technical: we construct an endomorphism $P$
of $V^{\o k}\o W^{\o\ell}$ (for each $(k,\ell)\in\N^2$) whose properties
will be used in the next section.

In section~\ref{simplesection}, we finally show that if $\L$ is a
nonelliptic BQD, then $\A_\L$ is indeed a quantum $\SL(3)$ (where the
$V_{(k,\ell)}$ of definition~\ref{Gqdef} are those appearing in the
$\N^2$-graded algebra $\M_\L$).

In section~\ref{equivsection}, we define an equivalence relation for BQD's
and then summarize the previous results, showing that---away from the
elliptic case---the correspondences $\A\mapsto\L_\A$ and $\L\mapsto\A_\L$
between quantum $\SL(3)$'s (up to isomorphism) and BQD's (up to
equivalence) are inverse of each other.

This raises the question of classifying BQD's. A related classification
problem has been studied in~\cite{EO}, where necessary conditions are
considered for a quantum analogue of $\O(\GL(3))$ to have correct
dimensions in degrees $\le4$. It turns out that these conditions, plus a
quantum determinant being central, plus a parameter not being a root of
unity, amount to our definition of a BQD.

Since we need an explicit classification of BQD's for a crucial case by
case argument in section~\ref{shapesection}, we reproduce it here, in
section~\ref{classsection}. This classification is complete, except for
case I.h, related to elliptic curves. (This case is however shown to
exist.) By the results of section~\ref{equivsection}, this also yields a
classification of all (nonelliptic) quantum $\SL(3)$'s. An important
ingredient in this classification will be a $3\times3$-matrix $Q$, which
encodes the square of the antipode and which can take four different Jordan
normal forms. The first possible form is the identity; we give a geometric
description of some cases there, in terms of plane cubic curves. The second
possible form has three different eigenvalues; it leads in particular to
the Artin-Schelter-Tate quantum $\SL(3)$'s~\cite{AST} (of which the
standard quantum $\SL(3)$~\cite{FRT} is a special case), and to the
Cremmer-Gervais one (see~\cite{Ho}). The third and fourth forms are
nondiagonal.

Finally, we list some indications for further study in
section~\ref{outlooksection}.

\section{Notations and conventions}

We denote by $\Z$ (resp.~$\N$, $\C$) the set of integers (resp.\
nonnegative integers, complex numbers). If $n\in\N$, $n\ge1$, and $t\in\C$,
let
\begin{align*}
[n]_t&:=1+t+\dots+t^{n-1}\\
[n]_t\,!&:=[1]_t\,[2]_t\,\dots\,[n]_t
\end{align*}
All vector spaces, algebras and tensor products are over $\C$. If $X,Y$ are
finite-dimensional vector spaces, we denote by $\Lin(X,Y)$ the space of
linear maps from $X$ to $Y$, and if $\al\in\Lin(X,Y)$, we denote by
$\tp\al\in\Lin(Y^*,X^*)$ its transpose. The identity map on $X$ is denoted
by $1_X$ (or simply~$1$). The tensor algebra of $X$ is denoted by $TX$. The
tensor product of $X$ and $Y$ will be denoted as usual by $X\o Y$, but for
typographical reasons, we denote the tensor product of two linear maps
$\al,\be$ by $(\al,\be)$.

If $\A$ is an algebra, we think of an element $\al\in\Lin(X,Y)\o\A$ as a
``linear map with coefficients in $\A$'', and we call \emph{space of
coefficients} of $\al$ the unique minimal vector subspace $\Coeff(\al)$ of
$\A$ such that $\al\in\Lin(X,Y)\o\Coeff(\al)$; obviously,
$\dim\Coeff(\al)\le(\dim X)(\dim Y)$.

Note that an equality in $\Lin(X,Y)\o\A$ amounts to $(\dim X)(\dim Y)$
equalities in $\A$. We shall use this to write relations in $\A$ in a
condensed way.

If $\al\in\Lin(X,Y)\o\A$ and $\be\in\Lin(Y,Z)\o\A$, there is an obvious
notion of composite $\be\al\in\Lin(X,Z)\o\A$ (using the multiplication in
$\A$). Similarly, if $\al\in\Lin(X,Y)\o\A$ and $\be\in\Lin(X',Y')\o\A$,
there is an obvious notion of tensor product $(\al,\be)\in\Lin(X\o X',Y\o
Y')\o\A$ (ditto).

All Hopf algebras are supposed to have an invertible antipode. ``Comodule''
means ``finite-dimensional right comodule''.

Let $\A$ be a Hopf algebra with comultiplication $\Delta$, counit $\eps$
and antipode $S$. We view an $\A$-comodule structure on a
finite-dimensional vector space $X$ as an element $t\in\Lin(X,X)\o\A$ such
that $\Delta(t)=t\o t$ and $\eps(t)=1_X$ (these are equalities in
$\Lin(X,X)\o(\A\o\A)$ and in $\Lin(X,X)$, respectively). Recall that every
$\A$-comodule is a direct sum of simple ones if and only if $\A$ is the
(direct) sum of the coefficient spaces of all (equivalence classes of)
simple $\A$-comodules. If so, this direct sum is called the
\emph{Peter-Weyl decomposition} of $\A$.

An $\A$-comodule morphism (more simply called $\A$-morphism) between two
$\A$-comodules $(X,t)$ and $(Y,u)$ is just an (ordinary) linear map
$\al\in\Lin(X,Y)$ such that $\al t=u\al$, where composites are taken in the
above sense. Tensor products of comodules also coincide with tensor
products in the above sense. The left dual of a comodule $(X,t)$ is
$(X^*,{}^*t)$, where ${}^*t=S(\!\tp\,t)\in\Lin(X^*,X^*)\o\A$, and the right
dual is $(X^*,t^*)$, with $t^*=S\inv(\!\tp\,t)$. Recall that for the left
(resp.\ right) dual structure, the canonical maps $X^*\o X\to\C$ and $\C\to
X\o X^*$ (resp.\ $X\o X^*\to\C$ and $\C\to X^*\o X$) are $\A$-morphisms.

\section{{}From quantum $\SL(3)$'s to BQD's}
\label{AtoLsection}

Let $\A$ be a quantum $\SL(3)$ (see definition~\ref{Gqdef}) and write
$V:=V_{(1,0)}$, $W:=V_{(0,1)}$.

\begin{prop}
\label{qSL3givesBQD}
There are $\A$-morphisms
\begin{equation}
\label{basic}
\begin{aligned}
A&:V\o V\to W\\
B&:W\o W\to V\\
C&:W\o V\to\C\\
D&:V\o W\to\C
\end{aligned}
\qquad
\begin{aligned}
a&:W\to V\o V\\
b&:V\to W\o W\\
c&:\C\to V\o W\\
d&:\C\to W\o V
\end{aligned}
\end{equation}
unique up to scalars, a constant $q\ne0$, $q^2\ne-1$, unique up to
$q\exch-q$ and $q\exch q\inv$, and a unique $3$-rd root of unity $\om$,
such that
\begin{subequations}
\label{coh}
\begin{align}
(1_V,C)(c,1_V)&=1_V&(D,1_V)(1_V,d)&=1_V\\
Aa&=1_W&&\\
C(A,1_V)&=\om\,D(1_V,A)&(1_V,a)c&=\om\,(a,1_V)d\\
(1_V,D)(a,1_W)&=B&\om^2\,(1_W,A)(d,1_V)&=b\\
\om\,(C,1_V)(1_W,a)&=B&(A,1_W)(1_V,c)&=b\\
Dc&=\kappa\,1_\C&Cd&=\kappa\,1_\C
\end{align}
\begin{align}
(1_V,A)(a,1_V)(A,1_V)(1_V,a)&=\rho\,(1_{V\o W}+cD)\\
(A,1_V)(1_V,a)(1_V,A)(a_V,1)&=\rho\,(1_{W\o V}+dC)
\end{align}
\end{subequations}
where $\ka=q^{-2}+1+q^2$ and $\rho=(q+q\inv)^{-2}$.
\end{prop}
(Recall that we write $(\al,\be)$ for the tensor product of any two linear
maps $\al$, $\be$.) Note that $q^2$ and $\om$ are the two parameters
of~\cite{KW} mentioned in the introduction.
\begin{proof}
First, definition~\ref{Gqdef}(c) asks for the following decompositions:
\begin{equation}
\label{decomp}
\begin{aligned}
V\o V&\simeq W\oplus V_{(2,0)}\\
W\o V&\simeq\C\oplus V_{(1,1)}
\end{aligned}
\qquad
\begin{aligned}
V\o W&\simeq\C\oplus V_{(1,1)}\\
W\o W&\simeq V\oplus V_{(0,2)}
\end{aligned}
\end{equation}
This already implies the desired uniqueness of the maps \eqref{basic}.

Furthermore, \eqref{decomp} prevents $V$ from being its own dual; so the
(left \emph{and} right) dual of $V$ must be $W$, and vice-versa. Hence
there exist $\A$-morphisms $C,c,D,d$ as in \eqref{basic} satisfying
(\ref{coh}a).

Again by \eqref{decomp}, there are $\A$-morphisms $A,a$ as in \eqref{basic}
satisfying (\ref{coh}b). By definition~\ref{Gqdef}(c), $V^{\o3}$ contains
exactly one copy of the trivial $\A$-comodule $\C$, hence
\[
\begin{aligned}
C(A,1_V)&=\la\,D(1_V,A)\\
(1_V,a)c&=\mu\,(a,1_V)d
\end{aligned}
\qquad(\la,\mu\in\C)
\]
We then define
\begin{align*}
B&:=(1_V,D)(a,1_W)=\mu\,(C,1_V)(1_W,a)\\
b&:=\la\mu\,(1_W,A)(d,1_V)=\la^2\mu\,(A,1_W)(1_V,c)
\end{align*}
As a consequence
\begin{align*}
Bb &=\la\mu\,(1_V,D)(a,1_W)(1_W,A)(d,1_V)
=\la\mu\,(1_V,D)(1_V,1_V,A)(a,1_V,1_V)(d,1_V)\\
&=(C,1_V)(1_V,A,1_V)(1_V,a,1_V)(1_V,c)
=(C,1_V)(1_V,c)
=1_V\\
Bb &=\la^2\mu^2(C,1_V)(1_W,a)(A,1_W)(1_V,c)
=\la^2\mu^2(C,1_V)(A,1_V,1_V)(1_V,1_V,a)(1_V,c)\\
&=\la^3\mu^3(1_V,D)(1_V,A,1_V)(1_V,a,1_V)(d,1_V)
=\la^3\mu^3(1_V,D)(d,1_V)
=\la^3\mu^3\,1_V
\end{align*}
so there is a $3$-rd root of unity $\om$ such that $\la\mu=\om^2$.
Moreover, any rescaling of the maps \eqref{basic} that leaves the relations
(\ref{coh}ab) intact also leaves $\la\mu$ invariant, hence $\om$ is unique.
Actually, there is such a rescaling after which $\la=\mu=\om$. This yields
(\ref{coh}cde). Now
\[
Cd=C(A,1_V)(a,1_V)d=D(1_V,A)(1_V,a)c=Dc
\]
which implies (\ref{coh}f). Next, define
\begin{equation}
\label{FG}
\begin{aligned}
F&:=(A,1_V)(1_V,a):V\o W\to W\o V\\
G&:=(1_V,A)(a,1_V):W\o V\to V\o W
\end{aligned}
\end{equation}
In view of \eqref{decomp}, we must have
\begin{equation}
\label{FG-GF}
\begin{aligned}
GF=(1_V,A)(a,1_V)(A,1_V)(1_V,a)&=\rho\,1_{V\o W}+\sig\,cD\\
FG=(A,1_V)(1_V,a)(1_V,A)(a_V,1)&=\rho'\,1_{W\o V}+\sig'\,dC
\end{aligned}
\end{equation}
for some $\rho,\sig,\rho',\sig'\in\C$, which are unique because any
rescaling that leaves (\ref{coh}abc) intact also leaves \eqref{FG-GF}
intact. Now use
\[
Fc=\om\,d\qquad CF=\om\,D\qquad Gd=\om^2\,c\qquad DG=\om^2\,C
\]
to compare $(GF)(GF)=G(FG)F$ and $(FG)(FG)=F(GF)G$: we get
$\rho'(\rho-\rho')=0$ and $\rho(\rho'-\rho)=0$, so $\rho'=\rho$.
Multiplying \eqref{FG-GF} on the right by $c$ (resp.~$d$) then yields
$1=\rho+\ka\sig$ and $1=\rho+\ka\sig'$, so $\sig'=\sig$. Next,
\begin{align*}
B(A,1_W)(1_V,G)&=\om^2\rho(D,1_V)+\om^2\sig(1_V,C)\\
=B(1_W,A)(F,1_V)&=\om^2\rho(1_V,C)+\om^2\sig(D,1_V)
\end{align*}
so $\sig=\rho$, which gives (\ref{coh}gh). Since $\rho\ne0$, the condition
$\rho=(q+q\inv)^{-2}$ defines $q$ with the desired uniqueness, and
$\ka=\rho\inv-1=q^{-2}+1+q^2$.
\end{proof}

\begin{prop}
\label{qnotroot}
In the notations of proposition~\ref{qSL3givesBQD}, either $q^2=1$, or
$q^2$ is not a root of unity.
\end{prop}

\begin{proof}
Let $R:=q\,1_{V\o V}-(q+q\inv)\,aA$. It follows from \eqref{coh} that
\begin{subequations}
\label{Hecke}
\begin{gather}
(R,1_V)(1_V,R)(R,1_V)=(1_V,R)(R,1_V)(1_V,R)\\
(R-q)(R+q\inv)=0
\end{gather}
\end{subequations}
If $k\ge2$, define the following endomorphisms of $V^{\o k}$:
\[
R_i:=1_{V^{\o(k-i-1)}}\o R\o 1_{V^{\o(i-1)}}\qquad(1\le i\le k-1)
\]
(This unusual right-to-left numbering will be convenient in
section~\ref{Heckesection}.)

Next, let us recall some general folklore on Hecke calculus. Let $\Sym_k$
be the symmetric group on $\{1,\dots,k\}$ and denote the transposition
$(i,i+1)$ by $s_i$ ($1\le i\le k-1$). If $w\in\Sym_k$ and if
$w=s_{i_1}\dots s_{i_p}$ is an expression of minimal length $p=:\ell(w)$,
let $R_w=R_{i_1}\dots R_{i_p}$; it follows from (\ref{Hecke}a) that $R_w$
does not depend on the choice of a particular such expression. Now let
\begin{equation}
\label{qsymm}
S_k:=\sum_{w\in\Sym_k}q^{\ell(w)}\,R_w
\end{equation}
It is easy to see that
\begin{equation}
\label{absorb}
R_iS_k=S_kR_i=q\,S_k
\end{equation}
for every $i$ (use the fact that $\ell(s_iw)=\ell(ws_i)=\ell(w)\pm1$ for
all $w$ and split the sum defining $S_k$ into two sums accordingly; then
use (\ref{Hecke}b)).

Recall also that every $w\in\Sym_k$ has a unique expression of minimal
length of the form $w=v_1v_2\dots v_{k-1}$, where
\begin{align*}
v_1&\in\{\,1\,,\,s_1\,\}\\
v_2&\in\{\,1\,,\,s_2\,,\,s_2s_1\,\}\\
&\;\vdots\\
v_{k-1}&\in\{1\,,\,s_{k-1},s_{k-1}s_{k-2}\,,\,\dots\,,\,s_{k-1}s_{k-2}\dots
s_1\,\}
\end{align*}
(This is just the bubble-sort principle.) It follows that
\begin{equation}
\label{Murphy}
\begin{aligned}
{}&S_k=(1+qR_1)(1+qR_2+q^2R_2R_1)\times\cdots\\
{}&\qquad\cdots\times(1+qR_{k-1}+\dots+q^{k-1}R_{k-1}\dots R_1)
\end{aligned}
\end{equation}
In particular, $S_k$ is of the form
\begin{equation}
\label{Sdecomp}
S_k=[k]_{q^2}\,!+a_{k-1}*+\dots+a_1*
\end{equation}
(where $a_i$ is defined the same way as $R_i$).

Now assume that $q^2$ is a primitive $n$-th root of unity, $n\ge2$. {}From
proposition~\ref{qSL3givesBQD}, we already know that $n\ge3$.

{\bf Claim~A.} \emph{If $2\le k\le n$, then $S_k\ne0$.}

\emph{First case: $k<n$.} As an $\A$-comodule, $\Im a_i$ is isomorphic to
$V^{\o(k-i-1)}\o W\o V^{\o(i-1)}$, so by definition~\ref{Gqdef}(c), $\Im
a_i$ does not contain $V_{(k,0)}$. It follows that $\sum_{i=1}^{k-1}(\Im
a_i)$ is a proper submodule of $V^{\o k}$. Now apply \eqref{Sdecomp},
noting that $[k]_{q^2}\,!\ne0$.

\emph{Second case: $k=n$.} (Adapted from~\cite[section~4]{KW}.) It follows
from \eqref{coh} that
\[
(1_{V^{\o(n-1)}},D)(R_{n-1},1_W)(1_{V^{\o(n-1)}},c)=q^3\,1_{V^{\o(n-1)}}
\]
Using this, \eqref{Murphy} and \eqref{absorb}, it follows that
\begin{align*}
(1_{V^{\o(n-1)}},D)(S_n,1_W)(1_{V^{\o(n-1)}},c)
&=(q^{-2}+1+q^2+q^4+\dots+q^{2n})\,S_{n-1}\\
&=(q^{-2}+1)\,S_{n-1}
\end{align*}
Since $S_{n-1}\ne0$ and $q^2\ne-1$, this shows claim~A.

By definition~\ref{Gqdef}(c), $V^{\o k}$ contains exactly one copy of
$V_{(k,0)}$; let us denote it by $M_{(k)}$.

{\bf Claim~B.} \emph{For $k\le n$, $\Im S_k=M_{(k)}$.}

(Proof adapted from~\cite[lemma~4.4]{KW}.) We proceed by induction over
$k$. If $k=2$, the statement is clear from the definitions, so assume
$2<k\le n$. Let $S'_k:=(S_{k-1},1_V)$ and $S''_k:=(1_V,S_{k-1})$. By
\eqref{Murphy}, $S_k=S'_kT'_k$ for some $T'_k$, and similarly,
$S_k=S''_kT''_k$, so $\Im S_k\subset\Im S'_k\cap\Im S''_k$. By induction,
$\Im S'_k=M_{(k-1)}\o V$ and $\Im S''_k=V\o M_{(k-1)}$. By
definition~\ref{Gqdef}(c), both images are isomorphic to the comodule
$V_{(k,0)}\oplus V_{(k-2,1)}$, so since $V^{\o k}$ contains only one copy
of $V_{(k,0)}$ (namely $M_{(k)}$), either $\Im S'_k\cap\Im S''_k=M_{(k)}$,
or $\Im S'_k=\Im S''_k$. But the second possibility would imply the
following equalities of subspaces in $V^{\o(2k)}$:
\[
M_{(k)}\o V^{\o k}=V\o M_{(k)}\o V^{\o(k-1)}=\dots=V^{\o k}\o M_{(k)}
\]
which is absurd. Therefore, $\Im S_k\subset M_{(k)}$, so by claim~A and by
the simplicity of $M_{(k)}$, this shows claim~B.

Now let $N_{(k)}$ be the unique supplemental comodule of $M_{(k)}$ in
$V^{\o k}$. {}From the proof of claim~A, it follows that
$\Im(S_k-[k]_{q^2}\,!)\subset N_{(k)}$ for all $k$. But for $k=n$, this
contradicts claim~B, because $[n]_{q^2}\,!=0$.
\end{proof}

\begin{df}
A \emph{basic quantum $\SL(3)$ datum} (BQD for short) consists of two
$3$-dimensional vector spaces $V$ and $W$, together with eight linear maps
\eqref{basic} satisfying \eqref{coh}, and such that either $q^2=1$, or
$q^2$ is not a root of unity.
\end{df}

If $\A$ is a quantum $\SL(3)$, we denote by $\L_\A$ the associated BQD.

A straightforward computation shows that in a BQD, relations (\ref{coh})
with $V\exch W$, $A\exch B$, $a\exch b$, $C\exch D$, $c\exch d$
interchanged are also satisfied, i.e.,
\begin{equation}
\label{postcoh}
\begin{gathered}
\begin{aligned}
(1_W,D)(d,1_W)&=1_W\\
Bb&=1_V\\
D(B,1_W)&=\om\,C(1_W,B)\\
(1_W,C)(b,1_V)&=A\\
\om\,(D,1_W)(1_V,b)&=A
\end{aligned}
\qquad
\begin{aligned}
(C,1_W)(1_W,c)&=1_W\\
&\\
(1_W,b)d&=\om\,(b,1_W)c\\
\om^2\,(1_V,B)(c,1_W)&=a\\
(B,1_V)(1_W,d)&=a
\end{aligned}\\
\begin{aligned}
(1_W,B)(b,1_W)(B,1_W)(1_W,b)&=\rho\,(1_{W\o V}+dC)\\
(B,1_W)(1_W,b)(1_W,B)(b,1_W)&=\rho\,(1_{V\o W}+cD)
\end{aligned}
\end{gathered}
\end{equation}

If $t\in\Lin(V,V)\o\A$ and $u\in\Lin(W,W)\o\A$ denote the $\A$-comodule
structures on $V$ and $W$, it follows from definition~\ref{Gqdef}(c) that
the coefficients of $t$ and $u$ generate $\A$. Moreover, the $\A$-morphisms
\eqref{basic} induce the following $4\times27+4\times9$ relations in $\A$:
\begin{equation}
\label{rels}
\begin{aligned}
A(t,t)&=uA\\
(t,t)a&=au
\end{aligned}
\qquad
\begin{aligned}
B(u,u)&=tB\\
(u,u)b&=bt
\end{aligned}
\qquad
\begin{aligned}
C(u,t)&=C\\
(t,u)c&=c
\end{aligned}
\qquad
\begin{aligned}
D(t,u)&=D\\
(u,t)d&=d
\end{aligned}
\end{equation}
(In coordinates, this reads $A^\al_{ij}t^i_kt^j_\ell=u^\al_\be
A^\be_{k\ell}$, etc.)

\section{{}From BQD's to Hopf algebras}
\label{LtoAsection}

Let us now work the other way round: if $\L$ is a BQD, the usual Tannakian
reconstruction procedure associates to it a bialgebra $\A_\L$, uniquely up
to unique isomorphism, together with $\A_\L$-comodule structures
$t\in\Lin(V,V)\o\A_\L$ and $u\in\Lin(W,W)\o\A_\L$, satisfying the two
following properties:
\begin{subequations}
\label{Tannaka}
\begin{align}
&\text{$A,a,B,b,C,c,D,d$ are $\A_\L$-morphisms}\\
&\text{$(\A,t,u)$ is universal with respect to (\ref{Tannaka}a)}
\end{align}
\end{subequations}
Condition (\ref{Tannaka}b) means that if $(\A',t',u')$ also satisfies
(\ref{Tannaka}a), then there is a unique bialgebra homomorphism
$\phi:\A_\L\to\A'$ such that $\phi(t)=t'$ and $\phi(u)=u'$.

Explicitely, $\A_\L$ is generated by the $(9+9)$-dimensional space
$\Coeff(t)+\Coeff(u)$, and relations \eqref{rels} form a presentation of
$\A_\L$.

The assignments $\Delta(t)=t\o t$, $\Delta(u)=u\o u$, $\eps(t)=1_V$,
$\eps(u)=1_W$ can be uniquely extended to algebra homomorphisms
$\D:\A_\L\to\A_\L\o\A_\L$ and $\eps:\A_\L\to\C$, which turn $\A_\L$ into a
bialgebra.

Moreover, relations (\ref{coh}a) turn $V$ and $W$ into each other's left
dual (in the monoidal category of $\A_\L$-comodules). This suggests an
antipode on $\A_\L$ defined by
\begin{equation}
\label{antipode}
S(t)=c^\flat\tp u\,C^\sharp\qquad S(u)=d^\flat\tp t\,D^\sharp
\end{equation}
(where $c^\flat$ is just $c$ viewed as a linear map $W^*\to V$, and
similarly for $C^\sharp:V\to W^*$, $d^\flat:V^*\to W$, $D^\sharp:W\to
V^*$). Indeed, using \eqref{coh} and \eqref{postcoh}, one checks that
\eqref{antipode} uniquely extends to an algebra antihomomorphism
$S:\A_\L\to\A_\L$, and this turns $\A_\L$ into a Hopf algebra.

\section{The shape algebra of a BQD}
\label{shapesection}

To a BQD $\L$, we associate the \emph{shape algebra} $\M_\L$, defined as
follows: let $I^G$ be the ideal in $T(V\oplus W)$ generated by all elements
$w\o v+(q+q\inv)\,G(w\o v)$, where $v\in V$, $w\in W$ (the map $G$ is
defined in \eqref{FG}). Then define $\M_\L:=T(V\oplus W)/I$, where $I$ is
the ideal generated by $\Im a$, $\Im b$, $\Im c$ and $I^G$. (Note that $I$
also contains $\Im d$ and all elements $v\o w+(q+q\inv)\,F(v\o w)$.)

The natural $\N^2$-grading on $\T(V\oplus W)$, with $V$ living in degree
$(1,0)$ and $W$ in degree $(0,1)$, factors to an $\N^2$-grading
$\M_\L=\bigoplus V_{(k,\ell)}$. There are natural identifications
$V_{(0,0)}\simeq\C$, $V_{(1,0)}\simeq V$, $V_{(0,1)}\simeq W$,
$V_{(2,0)}\simeq\Ker A$, $V_{(1,1)}\simeq\Ker C$, $V_{(1,1)}\simeq\Ker D$,
$V_{(0,2)}\simeq\Ker B$.

\begin{lem}
\label{Iinterright}
We have $\M_\L\simeq (TV\o TW)/I'$ (as $\N^2$-graded vector spaces), where
\[
I'=TV\o\Im a\o TV\o TW+TV\o\Im c\o TW+TV\o TW\o\Im b\o TW
\]
\end{lem}

\begin{proof}
First, the relations in $I^G$ can be turned into a reduction system (in the
sense of~\cite{Be}), which has no ambiguities at all; so by the diamond
lemma~\cite[thm~1.2]{Be},
\begin{equation}
\label{noamb}
T(V\oplus W)=(TV\o TW)\oplus I^G
\end{equation}
Next, if $w,w'\in W$, it follows from \eqref{coh} that
\[
(q+q\inv)^2\,(1,G)(G,1)\bigl(w\o a(w')\bigr)=a(w)\o w'+B(w\o w')\o c(1)
\]
This implies that
\begin{align*}
w\o a(w')&=a(w)\o w'+B(w\o w')\o c(1)\\
&\quad-(q+q\inv)\,\Bigl((q+q\inv)\,(1,G)(G,1)\bigl(w\o a(w')\bigr)
+(G,1)\bigl(w\o a(w')\bigr)\Bigr)\\
&\quad+\Bigl((q+q\inv)\,(G,1)\bigl(w\o a(w')\bigr)
+w\o a(w')\Bigr)
\end{align*}
Therefore
\[
W\o\Im a\subset\Im a\o W+V\o\Im c+I^G
\]
Similarly,
\begin{align*}
\Im b\o V&\subset V\o\Im b+\Im d\o W+I^G\\
\Im c\o V&\subset \Im a\o V+I^G\\
\Im d&\subset \Im c+I^G
\end{align*}
Applied inductively, these rules show that $I=I'+I^G$, hence $I\cap(TV\o
TW)=I'$ by \eqref{noamb}. The result follows.
\end{proof}

(There is of course a similar identification $\M_\L\simeq(TW\o TV)/I''$.)

\begin{lem}
\label{I3rightdim}
We have $\dim I_{(3,0)}=\dim I_{(0,3)}=17$ and $\dim I_{(2,1)}=\dim
I_{(1,2)}=66$.
\end{lem}

\begin{proof}
(We only look at $I_{(3,0)}$ and $I_{(2,1)}$.) Consider an element of $\Im
a\o V\cap V\o\Im a$, say $(a,1)(x)=(1,a)(y)$, where $x\in W\o V$ and $y\in
V\o W$. Applying $(A,1)$ and $(1,A)$ to this equality yields $x=F(y)$ and
$G(x)=y$, respectively, so $y=GF(y)=\rho\,y+\rho\,cD(y)$. But
$q^4+q^2+1\ne0$ implies $\rho\ne1$, hence $(1,a)(y)\in\Im(1,a)c$.
Conversely (see (\ref{coh}c)), $\Im(1,a)c=\Im(a,1)d\subset \Im a\o V\cap
V\o\Im a$. It follows that $I_{(3,0)}=\Im a\o V+V\o\Im a$ is of dimension
$9+9-1=17$.

Similarly, consider an element of $\Im a\o W\cap V\o\Im c$, say
$(a,1)(x)=(1,c)(y)$, where $x\in W\o W$ and $y\in V$. Applying $(A,1)$ and
$(1,D)$ yields $x=b(y)$ and $B(x)=\ka\,y$, so $y=Bb(y)=\ka\,y$. But
$q^4\ne-1$ implies $\ka\ne1$, so $y=0$. Thus, $I'_{(2,1)}=\Im a\o W+V\o\Im
c$ is of dimension $9+3-0=12$. Now apply lemma~\ref{Iinterright}.
\end{proof}

\begin{prop}
\label{Mright}
If $\L$ is not elliptic (i.e., not case I.h in section~\ref{classsection}),
then there are bases $(x_1,x_2,x_3)$ of $V$ and $(y_1,y_2,y_3)$ of $W$ such
that the monomials $x_1^ax_2^bx_3^cy_2^\ell y_1^k$ and
$x_1^ax_2^by_3^my_2^\ell y_1^k$ ($a,b,c,k,\ell,m\in\N$) form a basis of
$\M_\L$. Moreover, $\M_\L$ is a Koszul algebra.
\end{prop}
\begin{proof}
Using the classification given in section~\ref{classsection}, one can check
case by case that the relations coming from $\Im a$ may be written in the
form
\begin{align*}
x_3x_2&=\text{terms in $x_1^2$, $x_1x_2$, $x_1x_3$, $x_2x_1$, $x_2^2$,
$x_2x_3$, $x_3x_1$}\\
x_3x_1&=\text{terms in $x_1^2$, $x_1x_2$, $x_1x_3$, $x_2x_1$, $x_2^2$}\\
x_2x_1&=\text{terms in $x_1^2$, $x_1x_2$}
\intertext{those coming from $\Im b$ in the form}
y_2y_3&=\text{terms in $y_3^2$, $y_3y_2$}\\
y_1y_3&=\text{terms in $y_3^2$, $y_3y_2$, $y_3y_1$, $y_2y_3$, $y_2^2$}\\
y_1y_2&=\text{terms in $y_3^2$, $y_3y_2$, $y_3y_1$, $y_2y_3$, $y_2^2$,
$y_2y_1$, $y_1y_3$}
\intertext{that coming from $\Im c$ in the form}
x_3y_3&=\text{terms in $x_1y_1$, $x_1y_2$, $x_1y_3$, $x_2y_2$, $x_2y_3$}
\intertext{and those coming from $I^G$ in the form}
y_\al x_i&=\text{terms in $x_jy_\be$ ($j=1,2,3$, $\be=1,2,3$)}
\end{align*}
Now order the generators as follows: $x_1<x_2<x_3<y_3<y_2<y_1$; then order
the monomials degree-lexicographically. The relations defining $\M_\L$, as
written above, take the form of a reduction system (in the sense
of~\cite{Be}) that is compatible with the ordering just defined (i.e., each
term in a R.H.S. is strictly smaller than the corresponding L.H.S.). Since
the relations are homogeneous of degree~$2$, all ambiguities live in
degree~$3$, and since there are $50$ irreducible monomials of degree~$3$,
we know \emph{in advance} that all ambiguities are resolvable, thanks to
lemma~\ref{I3rightdim}. (This saves us from dozens of pages of ambiguity
cumputations!) Thus, the first statement follows from the diamond
lemma~\cite[thm~1.2]{Be}, noting that the irreducible monomials are exactly
those in the statement.

Finally, the basis so obtained is a labeled basis in the sense
of~\cite{Pr}, hence the last statement follows from~\cite[thm~5.3]{Pr}.
\end{proof}

The problem with the elliptic case (I.h in section~\ref{classsection}) is
that already the first three relations, viz.
\begin{align*}
\al x_2x_3+\be x_3x_2+\ga x_1^2&=0\\
\al x_3x_1+\be x_1x_3+\ga x_2^2&=0\\
\al x_1x_2+\be x_2x_1+\ga x_3^2&=0
\end{align*}
do not seem to be compatible with any semigroup ordering on the monomials,
so the diamond lemma does not apply.

\begin{cor}
\label{Mrightdim}
If $\L$ is not elliptic, $\dim
V_{(k,\ell)}=d_{(k,\ell)}=(k+1)(\ell+1)(k+\ell+2)/2$.
\end{cor}

Similarly, define the \emph{dual shape algebra} $\MM_\L:=\T(V^*\oplus
W^*)/J$, where $J$ is the ideal generated by $\Im\tp A$, $\Im\tp B$,
$\Im\tp D$ and $J^F$, the latter being the ideal generated by all elements
$\eta\o\xi+(q+q\inv)\,\tp F(\eta\o\xi)$ ($\xi\in V^*$, $\eta\in W^*$).

Again, we have an obvious $\N^2$-grading $\MM_\L=\bigoplus V^{(k,\ell)}$,
and the results of this section hold for $\MM_\L$ just as they do for
$\M_\L$.

{\bf Question.} When $\L$ is elliptic, is it still true that $\dim
V_{(k,\ell)}=\dim V^{(k,\ell)}=d_{(k,\ell)}$ and that $\M_\L$, $\MM_\L$ are
Koszul?

\section{Filtration of $\A_\L$ and dimensions}
\label{filtersection}

Consider the following subalgebra of $\MM_\L\o\M_\L$:
\[
\G_\L:=\bigoplus_{(k,\ell)\in\N^2} V^{(k,\ell)}\o V_{(k,\ell)}
\]
It is $\N^2$-graded (or $\N$-graded) by putting $V^{(k,\ell)}\o
V_{(k,\ell)}$ into degree~$(k,\ell)$ (or $k+\ell$).

\begin{prop}
The algebra $\G_\L$ is generated by $(V^*\o V)\oplus(W^*\o W)$, and
presented by the relations
\begin{equation}
\label{grArels}
\begin{gathered}
\begin{aligned}
A(t,t)&=0\\
(t,t)a&=0
\end{aligned}
\qquad
\begin{aligned}
B(u,u)&=0\\
(u,u)b&=0
\end{aligned}
\qquad
\begin{aligned}
D(t,u)&=0\\
(t,u)c&=0
\end{aligned}\\
(u,t)F-F(t,u)=0
\end{gathered}
\end{equation}
\end{prop}
(Here, $t$ denotes the canonical $(V^*\o V)$-comodule structure on $V$,
i.e., in coordinates, $t^i_j=x^i\o x_j$; similarly for $u$ w.r.t.\ $W$.)

\begin{proof}
Let $K$ be the ideal of $T\bigl((V^*\o V)\oplus(W^*\o W)\bigr)$ generated
by the L.H.S.'s of \eqref{grArels}: it follows from the definition of
$\M_\L$ and of $\MM_\L$ that we have an $\N^2$-graded map
\begin{equation}
\label{surjpresgrA}
T\bigl((V^*\o V)\oplus(W^*\o W)\bigr)/K\to\G_\L
\end{equation}
It already follows from lemma~\ref{Iinterright} that this map is
surjective. Since $F$ is invertible, we have
\[
T\bigl((V^*\o V)\oplus(W^*\o W)\bigr)=T(V^*\o V)\o T(W^*\o W)+K
\]
so for each $(k,\ell)\in\N^2$, \eqref{surjpresgrA} restricts to a
surjection
\begin{multline*}
{V^*}^{\o k}\o V^{\o k}\o{W^*}^{\o\ell}\o W^{\o\ell}
\Big/K\cap\bigl({V^*}^{\o k}\o V^{\o k}\o{W^*}^{\o\ell}\o
W^{\o\ell}\bigr)\\
\to
{V^*}^{\o k}\o{W^*}^{\o\ell}\o V^{\o k}\o W^{\o\ell}
\Big/\bigl({V^*}^{\o k}\o{W^*}^{\o\ell}\o I'_{(k,\ell)}
+J'_{(k,\ell)}\o V^{\o k}\o W^{\o\ell}\bigr)
\end{multline*}
(where $J'$ is to $\MM_\L$ what $I'$ is to $\M_\L$). Actually, this
surjection is just the obvious map induced from the flip $V^{\o
k}\o{W^*}^{\o\ell}\to{W^*}^{\o\ell}\o V^{\o k}$, so in view of the
definitions of $I'$, $J'$ and $K$, it is also injective.
\end{proof}

Note that~\eqref{grArels} is just the ``homogeneous part'' of the following
presentation of $\A_\L$:
\begin{equation}
\label{altrels}
\begin{gathered}
\begin{aligned}
A(t,t)&=uA\\
(t,t)a&=au
\end{aligned}
\qquad
\begin{aligned}
B(u,u)&=tB\\
(u,u)b&=bt
\end{aligned}
\qquad
\begin{aligned}
D(t,u)&=D\\
(t,u)c&=c
\end{aligned}\\
(u,t)F-F(t,u)=0
\end{gathered}
\end{equation}
so there is a canonical $\N$-graded surjection $\G_\L\to\gr\A_\L$, where we
filter $\A_\L=\bigcup_{n\ge0}\A^{(n)}$ by putting $\Coeff(t)+\Coeff(u)$
into degree~$1$.

\begin{prop}
\label{GisgrA}
If $\L$ is not elliptic, the surjection $\G_\L\to\gr\A_\L$ is an
isomorphism.
\end{prop}
\begin{proof}
Let $X:=(V^*\o V)\oplus(W^*\o W)$ and let $K_2$ be the subspace of $X\o X$
generated by the L.H.S.'s of~\eqref{altrels}, so $\G_\L=\T(X)/(K_2)$.
Consider the map $\la:K_2\to X$ sending each L.H.S. of~\eqref{altrels} to
the degree $1$ part of its R.H.S., and similarly $\mu:K_2\to\C$, for the
degree $0$ part; they are easily seen to be well-defined. Thanks
to~\cite[thm~0.5 and lemma~3.3]{BG}, we only have to prove the following
four conditions:
\begin{enumerate}
\item $\G_\L$ is Koszul,
\item the image of $\la\o1_X-1_X\o\la$ (defined on $K_2\o X\cap X\o K_2$)
lies in $K_2$,
\item $\la(\la\o1_X-1_X\o\la)=-(\mu\o1_X-1_X\o\mu)$,
\item $\mu(\la\o1_X-1_X\o\la)=0$.
\end{enumerate}
To prove condition (a), let us temporarily change signs in the grading of
$\MM_\L$, i.e., put $V^{(k,\ell)}$ into degree $(-k,-\ell)$. For the
resulting total $\Z^2$-grading on $\MM_\L\o\M_\L$, we now have
$\G_\L=(\MM_\L\o\M_\L)_{(0,0)}$. Since $\M_\L$ and $\MM_\L$ are Koszul
(proposition~\ref{Mright}), $\MM_\L\o\M_\L$ is also Koszul
by~\cite[prop.~2.1]{Pr}, and since the quadratic relations defining
$\MM_\L\o\M_\L$ are homogeneous w.r.t.~our temporary $\Z^2$-grading,
$\G_\L$ is still Koszul.

We now abandon this temporary grading and consider $\G_\L$ as an
$\N^2$-graded algebra, as before. Conditions (b), (c), (d) may clearly be
checked separately in degrees $(3,0)$, $(2,1)$, $(1,2)$, $(0,3)$. We only
look at degrees $(3,0)$ and $(2,1)$, since the other two are similar.

To improve legibility, we introduce the following notation: if $L$ is a
linear map, we write $L_i$ for a tensor product $(1,\dots,1,L,1,\dots,1)$,
with $L$ in the $i$-th place.

An element $\xi\in(K_2\o X\cap X\o K_2)_{(3,0)}$ is equal to both sides of
an equality
\[
\tr\Bigl[\al A_1(t,t,t)+(t,t,t)a_1\be\Bigr]
=
\tr\Bigl[\al'A_2(t,t,t)+(t,t,t)a_2\be'\Bigr]
\]
in $TX$, where
\begin{align*}
\al:W\o V&\to V\o V\o V&
\al':V\o W&\to V\o V\o V\\
\be:V\o V\o V&\to W\o V&
\be':V\o V\o V&\to V\o W
\end{align*}
(Here, $\tr$ means trace of a linear endomorphism with coefficients in
$\A_\L$; e.g., in coordinates, $\tr\bigl[\al A_1(t,t,t)\bigr]=\al_{\al
k}^{\ell mn}A_{ij}^\al\,t_\ell^i t_m^j t_n^k$.) It follows that
\begin{equation}
\label{split30}
\al A_1+a_1\be=\al'A_2+a_2\be'
\end{equation}
Now $2A_1\text{\eqref{split30}}a_1 -A_1\text{\eqref{split30}}a_2G
-FA_2\text{\eqref{split30}}a_1$ reads
\begin{align*}
&2(A_1\al+\be a_1)-FA_2\al-\be a_2G-A_1\al FG-FG\be a_1\\
&\quad=A_1\al'G+F\be'a_1-FA_2\al'G-F\be'a_2G
\end{align*}
Writing $-FG=FG-2FG=FG-2\rho-2\rho\,dC$, this becomes
\begin{align*}
2(1-\rho)(A_1\al+\be a_1)&=(A_1\al'+\be a_2-A_1\al F-F\be'a_2)G\\
&\quad+F(A_2\al+\be'a_1-A_2\al'G-G\be a_1)\\
&\quad+2\rho(A_1\al dC+dC\be a_1)
\end{align*}
Using an analogous expression of $2(1-\rho)(A_2\al'+\be'a_2)$, we get
\begin{align*}
&2(1-\rho)(\la\o1_X-1_X\o\la)(\xi)\\
&\quad=2(1-\rho)\tr\Bigl[(A_1\al+\be a_1)(u,t)\Bigr]
-2(1-\rho)\tr\Bigl[(A_2\al'+\be'a_2)(t,u)\Bigr]\\
&\quad=\tr\Bigl[(A_1\al'+\be a_2-A_1\al F-F\be'a_2)
\bigl(G(u,t)-(t,u)G\bigr)\Bigr]\\
&\qquad+\tr\Bigl[(A_2\al+\be'a_1-A_2\al'G-G\be a_1)
\bigl((u,t)F-F(t,u)\bigr)\Bigr]\\
&\qquad+2\rho\tr\Bigl[(A_1\al dC+dC\be a_1)(u,t)\Bigr]
-2\rho\tr\Bigl[(A_2\al'cD+cD\be'a_2)(t,u)\Bigr]
\end{align*}
Since the coefficients of $G(u,t)-(t,u)G$, $(u,t)F-F(t,u)$, $C(u,t)$,
$(u,t)d$, $D(t,u)$ and $(t,u)c$ are in $K_2$, this proves condition (b) in
degree $(3,0)$ (noting that $\rho\ne1$ follows from $q^4+q^2+1\ne0$).
Condition (c) is trivial. Next,
\begin{align*}
&2(1-\rho)\mu(\la\o1_X-1_X\o\la)(\xi)\\
&\quad=2\rho\tr\Bigl[A_1\al dC+dC\be a_1\Bigr]
-2\rho\tr\Bigl[A_2\al'cD+cD\be'a_2\Bigr]
\end{align*}
but the R.H.S. vanishes, as can be seen from multiplying \eqref{split30} on
the left by $CA_1=\om DA_2$ and on the right by $a_1d=\om^2a_2c$. This
proves condition (d) in degree $(3,0)$.

An element $\xi\in(K_2\o X\cap X\o K_2)_{(2,1)}$ is equal to both sides of
an equality
\begin{align*}
&\tr\Bigl[\al A_1(t,t,u)+(t,t,u)a_1\be\Bigr]
+\tr\Bigl[\ga C_1(u,t,t)+(u,t,t)d_1\de\Bigr]\\
&\qquad+\tr\Bigl[\eps(u,t,t)F_1-\eps F_1(t,u,t)\Bigr]\\
&\quad=\tr\Bigl[\al'A_2(u,t,t)+(u,t,t)a_2\be'\Bigr]
+\tr\Bigl[\ga'D_2(t,t,u)+(t,t,u)c_2\de'\Bigr]\\
&\qquad+\tr\Bigl[\eps'(t,t,u)G_2-\eps'G_2(t,u,t)\Bigr]
\end{align*}
where
\begin{align*}
\al:W\o W&\to V\o V\o W&
\al':W\o W&\to W\o V\o V\\
\be:V\o V\o W&\to W\o W&
\be':W\o V\o V&\to W\o W\\
\ga:V&\to W\o V\o V&
\ga':V&\to V\o V\o W\\
\de:W\o V\o V&\to V&
\de':V\o V\o W&\to V\\
\eps:W\o V\o V&\to V\o W\o V&
\eps':V\o V\o W&\to V\o W\o V
\end{align*}
It follows that
\begin{subequations}
\label{split21}
\begin{align}
\al A_1+a_1\be&=\ga'D_2+c_2\de'+G_2\eps'\\
\al'A_2+a_2\be'&=\ga C_1+d_1\de+F_1\eps\\
\eps F_1&=\eps'G_2
\end{align}
\end{subequations}
Noting that $G_1a_2=F_2a_1$, (\ref{split21}c)$G_1a_2$ reads
\[
\eps a_2+\om^2\eps d_1B=\eps'a_1+\eps'c_2B
\]
whereas $A_1$(\ref{split21}a)$a_1$ and $A_2$(\ref{split21}b)$a_2$ read
\begin{align*}
A_1\al+\be a_2&=A_1\ga'B+b\be'a_1+A_1G_2\eps'a_1\\
A_2\al'+\be'a_2&=\om^2A_2\ga B+\om b\de a_2+A_2F_1\eps a_2
\end{align*}
Since $A_1G_2=A_2F_1$, we therefore get
\begin{align*}
(\la\o1_X-1_X\o\la)(\xi)
&=\tr\Bigl[(A_1\al+\be a_1-A_2\al'-\be'a_2)(u,u)\Bigr]\\
&=\tr\Bigl[(A_1\ga'-\om^2A_2\ga+\om^2A_1G_2\eps
d_1-A_2F_1\eps'c_2)B(u,u)\Bigr]\\
&\quad+\tr\Bigl[(u,u)b(\de'a_1-\om\de a_2)\Bigr]
\end{align*}
Since the coefficients of $B(u,u)$ and $(u,u)b$ are in $K_2$, this proves
condition (b) in degree $(2,1)$. Next,
\begin{align*}
&\la(\la\o1_X-1_X\o\la)(\xi)\\
&\quad=\tr\Bigl[(BA_1\ga'-\om^2BA_2\ga+\om^2BA_1G_2\eps d_1-BA_2F_1\eps'c_2
+\de'a_1b-\de a_2b)\,t\Bigr]\\
&\quad=\tr\Bigl[\bigl(BA_1\ga'-\om^2BA_2\ga+\de'a_1b-\om\de
a_2b+\rho(D_1+C_2)
(\om\eps d_1-\om^2\eps'c_2)\bigr)\,t\Bigr]
\end{align*}
On the other hand,
\[
-(\mu\o1_X-1_X\o\mu)(\xi)
=\tr\Bigl[(D_2\ga'+\de'c_2-C_1\ga-\de d_1)\,t\Bigr]
\]
Now (\ref{split21}c)$c_1$ and (\ref{split21}c)$d_2$ read
\begin{align*}
\eps d_1&=\om^2\eps'G_2c_1=\om\eps'a_1b\\
\eps'c_2&=\om\eps F_1d_2=\eps a_2b
\end{align*}
so $(BA_1-D_2)$(\ref{split21}a)$(a_1b-c_2)$ reads
\[
0=(1-\ka)(BA_1\ga'-D_2\ga'+\de'a_1b-\de'c_2)
+\bigl(\rho D_1+(\rho-1)C_2\bigr)(\om\eps d_1-\om^2\eps'c_2)
\]
and $(\om^2BA_2-C_1)$(\ref{split21}b)$(\om a_2b-d_1)$ reads
\[
0=(1-\ka)(\om^2BA_2\ga-C_1\ga+\om\de a_2b-\de d_1)
+\bigl(\rho C_2+(\rho-1)D_1\bigr)(\om^2\eps'c_2-\om\eps d_1)
\]
Subtracting these two equalities and dividing by $1-\ka=(2\rho-1)/\rho$, we
get
\begin{align*}
0&=BA_1\ga'-D_2\ga'-\om^2BA_2\ga+C_1\ga+\de'a_1b-\de'c_2-\om\de a_2b+\de
d_1\\
&\quad+\rho(D_1+C_2)(\om\eps d_1-\om^2\eps'c_2)
\end{align*}
(Note that $\ka\ne1$ follows from $q^4\ne-1$.) This proves condition (c)
in degree $(2,1)$. Condition (d) is trivial.
\end{proof}

\begin{cor}
\label{AisPBW}
If $\L$ is not elliptic, $\dim\A^{(n)}=\sum_{k+\ell\le n}(d_{(k,\ell)})^2$.
\end{cor}
\begin{proof}
Use proposition~\ref{GisgrA} and corollary~\ref{Mrightdim} (applied to
$\M_\L$ and to $\MM_\L$).
\end{proof}

\section{A key endomorphism}
\label{Heckesection}

Let
\[
R:=q-(q+q\inv)\,aA\qquad\bR:=q\inv-(q\inv+q)\,bB
\]
(We have already used $R$ in the proof of proposition~\ref{qnotroot}.) Fix
$(k,\ell)\in\N^2$ and define the following endomorphisms of $V^{\o k}\o
W^{\o\ell}$:
\begin{align*}
R_i&:=1_{V^{\o(k-i-1)}}\o R\o 1_{V^{\o(i-1)}}\o 1_{W^{\o\ell}}&
&\qquad 1\le i\le k-1\\
X&:=1_{V^{\o(k-1)}}\o cD\o 1_{W^{\o(\ell-1)}}\\
\bR_i&:=1_{V^{\o k}}\o 1_{W^{\o(i-1)}}\o \bR \o 1_{W^{\o(\ell-i-1)}}&
&\qquad 1\le i\le \ell-1
\end{align*}

\begin{prop}
\label{keyendo}
The ring of endomorphisms of $V^{\o k}\o W^{\o\ell}$ generated by the
$R_i$'s, the $\bR_i$'s and $X$ contains an element $P$ such that the kernel
of the multiplication $V^{\o k}\o W^{\o\ell}\to V_{(k,\ell)}$ is contained
in $\Ker P$ and contains $\Im(P-1)$.
\end{prop}

\begin{proof}
First, the following relations are easily obtained from \eqref{coh} and
\eqref{postcoh}:
\begin{equation}
\label{RXRrels}
\begin{gathered}
\begin{aligned}
R_iR_j&=R_jR_i&
\quad\bR_i\bR_j&=\bR_j\bR_i&
\qquad&\text{if $|i-j|\ge2$}\\
XR_i&=R_iX&
X\bR_i&=\bR_iX&
&\text{if $i\ge2$}
\end{aligned}
\\[2mm]
\begin{aligned}
R_i\bR_j&=\bR_jR_i\\
R_iR_{i+1}R_i&=R_{i+1}R_iR_{i+1}\\
\bR_i\bR_{i+1}\bR_i&=\bR_{i+1}\bR_i\bR_{i+1}\\
(R_i-q)(R_i+q\inv)&=0\\
(\bR_i-q\inv)(\bR_i+q)&=0
\end{aligned}
\qquad
\begin{aligned}
X^2&=\kappa X\\
XR_1X&=q^3X\\
X\bR_1X&=q^{-3}X\\
XR_1\bR_1XR_1&=XR_1\bR_1X{\bR_1}\inv\\
\bR_1XR_1\bR_1X&=R_1\inv XR_1\bR_1X
\end{aligned}
\end{gathered}
\end{equation}
Note that these are the relations appearing in~\cite[def.~2.1]{KM}.
(\emph{Warning:} the index convention used here differs from~\cite{KM}, and
also from that in the proof of proposition~\ref{GisgrA}.) In the sequel, we
shall use the letters $a,b,c,d$ as indices; this should not cause confusion
with the maps \eqref{basic}.

Let $T^a_b:=\bR_a\dots\bR_1 XR_1\dots R_b$ (for $0\le a\le \ell-1$, $0\le
b\le k-1$). Using \eqref{RXRrels}, we get
\begin{equation}
\label{RTrels}
\begin{aligned}
T^a_b R_p&=
\begin{cases}
R_{p+1}T^a_b&\text{if $p\le b-1$}\\
(q-q\inv)T^a_b+T^a_{b-1}&\text{if $p=b$}\\
T^a_{b+1}&\text{if $p=b+1$}\\
R_pT^a_b&\text{if $p\ge b+2$}
\end{cases}
\\[2mm]
\bR_pT^a_b&=
\begin{cases}
T^a_b\bR_{p+1}&\text{if $p\le a-1$}\\
(q\inv-q)T^a_b+T^{a-1}_b&\text{if $p=a$}\\
T^{a+1}_b&\text{if $p=a+1$}\\
T^a_b\bR_p&\text{if $p\ge a+2$}
\end{cases}
\end{aligned}
\end{equation}
and
\begin{equation}
\label{TTrels}
\begin{aligned}
T^a_bT^c_d&=R_1\inv T^{c-1}_bT^a_d
&\qquad&\text{if $a\ge c\ge 1$ and $b\ge1$}\\
T^a_bT^c_d&=T^a_dT^c_{b-1}{\bR_1}\inv
&&\text{if $d\ge b\ge 1$ and $c\ge 1$}\\
T^a_bT^0_d&=q^3\,R_2\dots R_b\,T^a_d
&&\text{if $b\ge1$}\\
T^a_0T^c_d&=q^{-3}\,T^a_d\,\bR_c\dots\bR_2
&&\text{if $c\ge1$}\\
T^a_0T^0_d&=\kappa\,T^a_d\\
\end{aligned}
\end{equation}
Define $S=S_k$ as in \eqref{qsymm}, and define $\bS=\bS_\ell$ similarly
(replacing $R_i$ by $\bR_i$ and $q$ by $q\inv$). For $m\ge1$, let
\[
U_m:=\sum_{\substack{
0\le a_1<\dots<a_m\le\ell-1\\[1mm]
k-1\ge b_1>\dots>b_m\ge0
}}
q^{(b_1-a_1)+\dots+(b_m-a_m)}\,T^{a_1}_{b_1}\dots T^{a_m}_{b_m}
\]
(Note that $U_m=0$ if $m>\min(k,\ell)$.) Consider a linear combination
\[
P:=\al_0\,S\bS+\sum_{m=1}^{\min(k,\ell)}\al_m\,S U_m \bS
\]
Since $\bS\bR_i=q\inv\bS$, we have $P\bR_i=q\inv P$ ($1\le i\le\ell-1$). It
also follows from \eqref{RTrels} and \eqref{absorb} that
$SU_m\bS\,R_i=q\,SU_m\bS$ for $1\le i\le k-1$, so $PR_i=q\,P$.

Let $\bSp$ be $S^*_{\ell-1}$ acting on the $\ell-1$ rightmost factors of
$V^{\o k}\o W^{\o\ell}$; we still have $\bSp\bR_i=q\inv\bSp$ for $i\ge2$.
Let also $U'_m$ be the sum of all terms of $U_m$ in which $b_m=0$, and let
$U''_m=U_m-U'_m$. By an equality analogous to \eqref{Murphy}, we have
\begin{equation}
\label{Murphyswitch}
\bS X=(T^0_0+q\inv\,T^1_0+\dots+q^{-(\ell-1)}\,T^{\ell-1}_0)\bSp
\end{equation}
{}From \eqref{RTrels}, \eqref{TTrels} and \eqref{absorb}, it follows that
\begin{align*}
SU'_m(T^0_0+q\inv\,T^1_0+\dots+q^{-(\ell-1)}\,T^{\ell-1}_0)\bSp
&=q^2\,[\ell+2]_{q^{-2}}\,SU'_m\bSp\\
SU''_m(T^0_0+q\inv\,T^1_0+\dots+q^{-(m-1)}\,T^{m-1}_0)\bSp
&=q^4\,[k-m]_{q^2}\,SU'_m\bSp\\
SU''_m
(q^{-m}\,T^m_0+q^{-(m+1)}\,T^{m+1}_0+\dots+q^{-(\ell-1)}\,T^{\ell-1}_0)\bSp
&=[m+1]_{q^{-2}}\,SU'_{m+1}\bSp
\end{align*}
Adding these three relations and using \eqref{Murphyswitch}, we get
\[
SU_m\bS\,X=q^{-2\ell}\,[k+\ell-m+2]_{q^2}\,S U'_m \bSp
+[m+1]_{q^{-2}}\,S U'_{m+1} \bSp
\]
It follows that for a suitable choice of the constants $\al_m$, we have
$PX=0$. As a reminder, we have so far obtained the relations
\begin{equation}
\label{Pprop1}
PR_i=q\,P\qquad P\bR_i=q\inv\,P\qquad PX=0
\end{equation}
Next, if we define the maps $a_i$ and $b_i$ the same way as $R_i$ and
$\bR_i$, respectively, then $P$ is of the form
\begin{equation}
\label{Pprop2}
P=\al_0\,[k]_{q^2}\,!\,[\ell]_{q^{-2}}\,!\,
+a_{k-1}*+\dots+a_1*+c*+b_1*+\dots+b_{\ell-1}*
\end{equation}
Since $[n]_{q^2}\ne0$ for every $n$, we may rescale the $\al_m$ so that
$\al_0\,[k]_{q^2}\,!\,[\ell]_{q^{-2}}\,!=1$. Now combine \eqref{Pprop1},
\eqref{Pprop2} and lemma~\ref{Iinterright}.
\end{proof}

\section{Simple $\A_\L$-comodules}
\label{simplesection}

\begin{prop}
\label{allsimples}
If $\L$ is a nonelliptic BQD, then $\A_\L$ is a quantum $\SL(3)$.
\end{prop}

\begin{proof}
In the sequel, reasonings involving a degree $(k,\ell)\in\N^2$ will work
even in the limit cases $k=0$ and $\ell=0$, provided one drops anything
involving a negative index.

The canonical $\bigl((V^*\o V)\oplus(W^*\o W)\bigr)$-comodule structure on
$V\oplus W$ turns the shape algebra $\M_\L$ into an $\N^2$-graded
(infinite-dimensional) $\A_\L$-comodule algebra, hence every $V_{(k,\ell)}$
($(k,\ell)\in\N^2$) into an $\A_\L$-comodule. Note that by
corollary~\ref{Mrightdim}, the dimension requirements of
definition~\ref{Gqdef} are satisfied.

We first prove the following statements by induction over $n$:
\begin{itemize}
\item[A$_n$.] All $V_{(k,\ell)}$, $k+\ell\le n$, are simple and pairwise
nonisomorphic.
\item[B$_n$.] For every $(k,\ell)$, $k+\ell\le n$, we have decompositions
\begin{subequations}
\label{CG}
\begin{align}
V_{(k,\ell)}\o V
&\simeq V_{(k+1,\ell)}\oplus V_{(k,\ell-1)}\oplus V_{(k-1,\ell+1)}\\
V_{(k,\ell)}\o W
&\simeq V_{(k,\ell+1)}\oplus V_{(k-1,\ell)}\oplus V_{(k+1,\ell-1)}
\end{align}
\end{subequations}
\end{itemize}
Statements A$_0$~and B$_0$ are clear.

Let $n>0$. Applying statement~B$_{n-1}$ repeatedly, we get
\[
\bigoplus_{p=0}^n(V\oplus W)^{\o p}\simeq V_{(k_1,\ell_1)}\oplus\dots\oplus
V_{(k_s,\ell_s)}
\]
where $k_i+\ell_i\le n$ for every $i$ (and actually, every $V_{(k,\ell)}$,
$k+\ell\le n$, appears at least once), so the coefficient space of the
R.H.S. has dimension at most $\sum_{k+\ell\le n}(d_{(k,\ell)})^2$, with
equality if and only if statement~A$_n$ holds. But in view of the L.H.S.,
this coefficient space is precisely $\A^{(n)}$, so by
corollary~\ref{AisPBW}, we do have the desired equality. This shows that
B$_{n-1}$ implies A$_n$.

Let still $n>0$ and let $k+\ell=n$. Denote by $\mu$ the multiplication in
$\M_\L$ and define
\begin{align*}
\al&:V_{(k-1,\ell)} \xto{(1,c)}
V_{(k-1,\ell)}\o V\o W \xto{(\mu,1)}
V_{(k,\ell)}\o W\\
\be&:V_{(k+1,\ell-1)} \xto{\ga}
V_{(k,\ell-1)}\o V \xto{(1,b)}
V_{(k,\ell-1)}\o W\o W \xto{(\mu,1)}
V_{(k,\ell)}\o W
\end{align*}
where $\ga$ is the injection (provided by statement~B$_{n-1}$) such that
$\mu\ga=1$. Consider the sequence
\begin{equation}
\label{CGexact}
0 \to V_{(k-1,\ell)}\oplus V_{(k+1,\ell-1)}
\xto{\al\oplus\be}
V_{(k,\ell)}\o W
\xto\mu
V_{(k,\ell+1)} \to 0
\end{equation}
By lemma~\ref{Iinterright}, $\mu$ is surjective. Moreover,
\[
(1,C)(\al,1_V)=(1,C)(\mu,1,1)(1,c,1)=\mu(1,1,C)(1,c,1)=\mu
\]
and
\begin{align*}
\mu(1,B)(\be,1_W)&=\mu(1,B)(\mu,1,1)(1,b,1)(\ga,1)
=\mu(\mu,1)(1,1,B)(1,b,1)(\ga,1)\\
&=\om\,\mu(\mu,1)(1,F)(\ga,1)
=\om\,\mu(1,\mu)(1,F)(\ga,1)\\
&=\om(q+q\inv)\inv\,\mu(1,\mu)(\ga,1)
=\om(q+q\inv)\inv\,\mu(\mu,1)(\ga,1)\\
&=\om(q+q\inv)\inv\,\mu
\end{align*}
hence $\al\ne0$ and $\be\ne0$. Statement~A$_n$ then implies that
$\al\oplus\be$ is injective. Moreover,
\begin{align*}
\mu\al&=\mu(\mu,1)(1,c)=\mu(1,\mu)(1,c)=0\\
\mu\be&=\mu(\mu,1)(1,b)\ga=\mu(1,\mu)(1,b)\ga=0
\end{align*}
so in view of dimensions, the sequence \eqref{CGexact} is exact. Using the
$\A$-endomorphism $P$ provided by proposition~\ref{keyendo}, we may split
this exact sequence; this implies (\ref{CG}b). Since (\ref{CG}a) is
similar, this shows that B$_{n-1}$ and A$_n$ imply B$_n$ and completes the
induction.

Definition~\ref{Gqdef}(a) is now clear. Moreover, the sum of the
coefficient spaces of all the $V_{(k,\ell)}$ is equal to
$\bigcup\A^{(n)}=\A_\L$, which implies definition~\ref{Gqdef}(b).

Finally, (\ref{CG}a) implies
\begin{equation}
\label{CGinduct}
\begin{aligned}
{}&V_{(k,\ell)}\o V_{(p,q)}\o V\\
{}&\quad\simeq(V_{(k,\ell)}\o V_{(p+1,q)})
\oplus(V_{(k,\ell)}\o V_{(p,q-1)})
\oplus(V_{(k,\ell)}\o V_{(p-1,q+1)})
\end{aligned}
\end{equation}
By induction over $p+q$, $V_{(k,\ell)}\o V_{(p,q-1)}$, $V_{(k,\ell)}\o
V_{(p-1,q+1)}$ and $V_{(k,\ell)}\o V_{(p,q)}$ decompose according to
definition~\ref{Gqdef}(c), and by \eqref{CG}, $V_{(k,\ell)}\o V_{(p,q)}\o
V$ still decomposes according to it. Therefore, \eqref{CGinduct} implies
that $V_{(k,\ell)}\o V_{(p+1,q)}$ does so as well. The case $V_{(k,\ell)}\o
V_{(p,q+1)}$ is similar.
\end{proof}

\section{Equivalence between quantum $\SL(3)$'s and BQD's}
\label{equivsection}

Let $\L=(V,W,A,a,B,b,C,c,D,d)$ be a BQD and consider the following three
transformations of $\L$:
\begin{itemize}
\item base change, i.e., conjugating $A,a,B,b,C,c,D,d$ by some invertible
linear maps $V\to V'$, $W\to W'$ (where $V',W'$ are any vector spaces of
dimension~$3$);
\item multiplying $A,a,B,b,C,c,D,d$ by scalars (but such that \eqref{coh}
are preserved);
\item interchanging $V\exch W$, $A\exch B$, $a\exch b$, $C\exch D$, $c\exch
d$ (this still gives a BQD by \eqref{postcoh}).
\end{itemize}
We call the third transformation \emph{Dynkin flip}, and it should indeed
be thought of as applying the automorphism of the Dynkin diagram of
$\SL(3)$. Two BQD's are called \emph{equivalent} if one is obtained from
the other by any combination of these three transformations.

\begin{theo}
The correspondences $\A\mapsto\L_\A$ and $\L\mapsto\A_\L$ are inverse of
each other between nonelliptic quantum $\SL(3)$'s (up to Hopf algebra
isomorphism) and nonelliptic BQD's (up to equivalence).
\end{theo}

\begin{proof}
These two correspondences are well-defined by
propositions~\ref{qSL3givesBQD}, \ref{qnotroot}~and \ref{allsimples}. Also,
if $\L$ is not elliptic, then $\L_{\A_\L}$ is clearly equivalent to $\L$.

Conversely, if $\A$ is a nonelliptic quantum $\SL(3)$, then the algebras
$\A_{\L_\A}$ and $\A$ have a common generating space $\Coeff(t)+\Coeff(u)$,
so both have an $\N$-filtration $\A_{\L_\A}=\bigcup\A^{(n)}$,
$\A=\bigcup\B^{(n)}$. Furthermore, the defining relations \eqref{rels} of
$\A_{\L_\A}$ are also valid in $\A$, so there is a canonical surjection
$\A_{\L_\A}\to\A$, which restricts to the identity on
$\Coeff(t)+\Coeff(u)$, and therefore takes $\A^{(n)}$ onto $\B^{(n)}$. But
\[
\dim\A^{(n)}=\sum_{k+\ell\le n}(d_{(k,\ell)})^2=\dim\B^{(n)}
\]
where the left equality follows from corollary~\ref{AisPBW} and the right
one from the Peter-Weyl decomposition. Thus, $\A_{\L_\A}$ is isomorphic to
$\A$.
\end{proof}

\section{Classification of BQD's}
\label{classsection}
\setcounter{subsection}{-1}

\subsection{Strategy for the classification}

Let
\[
Q:=c^\flat D^\sharp:V\to V
\]
By (\ref{coh}a), $Q$ is invertible with $Q\inv=d^\flat C^\sharp$, so by
(\ref{coh}f),
\begin{equation}
\label{trc}
\tr Q=\ka=\tr Q\inv
\end{equation}
(Note that by \eqref{antipode}, $Q$ ``encodes'' the square of the antipode,
in the sense that $S^2(t)=QtQ\inv$.) Now define
\begin{equation}
\label{Edef}
E:=C(A,1_V)=\om\,D(1_V,A)
\qquad
e:=(1_V,a)c=\om\,(a,1_V)d
\end{equation}
Let $x_1,x_2,x_3$ be a basis of $V$ and $y_1,y_2,y_3$ a basis of $W$, and
denote by $x^1,x^2,x^3$ and $y^1,y^2,y^3$ their respective dual bases.
Write $Q(x_i)=Q^j_i x_j$, $e=e^{ijk}x_{ijk}$ and $E=E_{ijk}x^{ijk}$, where
\[
x_{ijk}:=x_i\o x_j\o x_k
\qquad
x^{ijk}:=x^i\o x^j\o x^k
\]
By \eqref{coh}, we have
\begin{gather}
\label{QE}
e^{\ell ij}=\om^2\,e^{ijk}Q^\ell_k
\qquad
E_{\ell ij}=\om\,E_{ijk}Q^k_\ell\\
\label{Ee}
e^{ik\ell}E_{k\ell j}=\delta^i_j
\end{gather}

\begin{prop}
$Q$ is of one of the following four types:
\begin{description}
\item[Type I] $Q$ is the identity.
\item[Type II] $Q$ has eigenvalues $q^{-2},1,q^2$ for some $q$, with
$q^{2n}\ne1$ for
all $n$ (so $Q$ is diagonalizable).
\item[Type III] $Q$ has triple eigenvalue~$1$ and a $2\times2$ Jordan
block.
\item[Type IV] $Q$ has triple eigenvalue~$1$ and a $3\times3$ Jordan block.
\end{description}
\end{prop}

\begin{proof}
Let $\al,\be,\ga$ be the three eigenvalues of $Q$.

\emph{Case~1:} $\al=\be=\ga$. Then (up to base change)
$Q=\spmat{\al&0&0\\0&\al&0\\0&0&\al}$,
$\spmat{\al&0&2\\0&\al&0\\0&0&\al}$ or
$\spmat{\al&2&2\\0&\al&2\\0&0&\al}$.
(The unusual normalizations of the nondiagonal cases will be convenient
later.)

By \eqref{trc}, $\al^2=1$. Taking $\al=-1$ in \eqref{QE} prevents $e,E$
from satisfying \eqref{Ee} (for either value of $\om$ and either form of
$Q$), so we must have $\al=1$. This gives types I, III~and IV.

\emph{Case~2:} $\al=\ga\ne\be$. Then $Q=\spmat{\al&0&0\\0&\be&0\\0&0&\al}$
or $\spmat{\al&0&1\\0&\be&0\\0&0&\al}$.

If $\al^2\be=1$, then \eqref{trc} implies $\al=\pm1$, hence $\be=1$, and
$\al=-1$. But then $\ka=-1$, contradicting $q^2\ne-1$. Therefore
$\al^2\be\ne1$. Taking $\al^3\ne1$ in \eqref{QE} then prevents $e,E$ from
satisfying \eqref{Ee}, so we must have $\al^3=1$.

If $\al\be^2=1$, then \eqref{trc} implies $\al=1$, $\be=-1$, hence $\ka=1$,
contradicting $q^4\ne-1$. Therefore $\al\be^2\ne1$.

Now \eqref{Ee} implies $e^{222}E_{222}=1$, so \eqref{QE} implies $\be=\om$
and $\be=\om^2$, hence $\be=1$. Then \eqref{trc} implies $\al^2=1$,
contradicting $\al^2\be\ne1$. Consequently, case~2 is impossible.

\emph{Case~3:} $\al\ne\be\ne\ga\ne\al$. Then
$Q=\spmat{\al&0&0\\0&\be&0\\0&0&\ga}$.

\emph{Subcase~3a:} $\al\be\ga=1$. Then \eqref{trc} implies
$(\al-1)(\be-1)(\ga-1)=0$, hence, say, $\be=1$. Then $\al=\ga\inv$ and
$\ka=\ga\inv+1+\ga$. Choosing $q$ such that $q^2=\ga$ (so that
$\ka=q^{-2}+1+q^2$) gives type~II.

\emph{Subcase~3b:} $\al\be\ga\ne1$. If, say, $\be=1$, then \eqref{trc}
implies $\ga=\al\inv$ or $\ga=-\al$. Both are impossible ($\ga=\al\inv$
contradicts $\al\be\ga\ne1$ and $\ga=-\al$ implies $\ka=1$, contradicting
$q^4\ne-1$), so we must have $\al\ne1$, $\be\ne1$ and $\ga\ne1$.

If at least five of the six scalars $\al^2\be$, $\al^2\ga$, $\be^2\al$,
$\be^2\ga$, $\ga^2\al$, $\ga^2\be$ are different from~$1$, then \eqref{QE}
prevents $e,E$ from satisfying \eqref{Ee}. Therefore, say, $\al^2\be=1$ (so
that $\al^2\ga\ne1$ and $\be^2\al\ne1$) and at least one of $\be^2\ga$,
$\ga^2\al$, $\ga^2\be$ equals~$1$.

But $\be^2\ga=1$ implies $\ga=\be^{-2}=\al^4$, so by \eqref{trc}, $\ka=0$,
contradicting $q^4+q^2+1\ne0$. Similarly, $\ga^2\al=1$ is impossible.

If $\ga^2\be=1$, then $\ga=-\al$, and \eqref{trc} implies $\al^4=1$, hence
$\be=\pm1$, so $\ka=\pm1$, both of which we already know to be excluded.
Consequently, subcase~3b is impossible.
\end{proof}
If $Q$ is of type X (X$=$I, II, III or IV), we call the BQD of type X if
$\om=1$, and of type X' if $\om^2+\om+1=0$.

The strategy to classify BQD's of type I, I', II, II', III, III', IV~and
IV' will be as follows:
\begin{itemize}
\item In types II and II', take the given $q$; in all other types, set
$q=1$.
\item The choice of a ``nondegenerate'' $c:\C\to V\o W$ being arbitrary
(changing it amounts to a base change in $W$), deduce $D$ from $Q=c^\flat
D^\sharp$, then $C,d$ from (\ref{coh}a).
\item Choose $e:\C\to V\o V\o V$ satisfying \eqref{QE}, working modulo the
stabilizer $\Stab(Q)$ of $Q$ in $\GL(V)$.
\item Determine all possible $E:V\o V\o V\to\C$ satisfying \eqref{QE},
\eqref{Ee} and (\ref{coh}gh) (where $A,a$ are deduced from \eqref{Edef}).
\item Reduce the possible forms for $E$ modulo $\Stab(Q,e)$.
\end{itemize}
In some cases, it may be useful to swap the last two steps.

We shall allow ourselves to satisfy \eqref{Ee} only up to a nonzero scalar,
adapting (\ref{coh}b) and (\ref{coh}gh) accordingly.

Finally, we leave out the details of matrix computations, which the reader
can easily recheck using any standard computer algebra package. (The author
used Maple~V Release~3.)

\subsection{Type~I}

We take
\[
\begin{aligned}
c&=x_1\o y_1+x_2\o y_2+x_3\o y_3\\
C&=y^1\o x^1+y^2\o x^2+y^3\o x^3
\end{aligned}
\qquad
\begin{aligned}
d&=y_1\o x_1+y_2\o x_2+y_3\o x_3\\
D&=x^1\o y^1+x^2\o y^2+x^3\o y^3
\end{aligned}
\]
We have $\Stab(Q)=\GL(V)$. By \eqref{QE}, $e=\la+s$ and $E=\La+S$, where
$\la,\La$ are totally antisymmetric and $s,S$ totally symmetric. We view
$s$ (resp.\ $S$) as a polynomial function on $V^*$ (resp.\ $V$), viz.
\begin{align*}
s&=e^{111}x_1^3+e^{222}x_2^3+e^{333}x_3^3+3(e^{123}+e^{132})x_1x_2x_3\\
&\quad+3e^{112}x_1^2x_2+3e^{113}x_1^2x_3
+3e^{221}x_2^2x_1+3e^{223}x_2^2x_3\\
&\quad+3e^{331}x_3^2x_1+3e^{332}x_3^2x_2\\
S&=E_{111}X_1^3+E_{222}X_2^3+E_{333}X_3^3+3(E_{123}+E_{132})X_1X_2X_3\\
&\quad+3E_{112}X_1^2X_2+3E_{113}X_1^2X_3
+3E_{221}X_2^2X_1+3E_{223}X_2^2X_3\\
&\quad+3E_{331}X_3^2X_1+3E_{332}X_3^2X_2
\end{align*}
where we write $X_i$ instead of $x^i$, to improve legibility. We shall also
consider $s$ as a cubic curve in the projective plane $\P V^*$ and $S$ as
one in $\P V$.

We now examine the different normal forms for $s$ modulo $\GL(V)$ if
$\la=0$, and modulo $\Stab(\la)=\SL(V)$ if $\la\ne0$ (see, e.g., \cite[\S
I.7]{Kr}).

The following cases are possible for $\la\ne0$, $\La\ne0$ (in cases
I.a--I.g, we normalize so that $\la_{123}=1$, $\La^{321}=1$).
\begin{description}
\item[Case I.a] $s=0$, $S=0$. The resulting Hopf algebra is that of
functions
on the
(ordinary) group $\SL(3)$.
\item[Case I.b] $s=x_1^3$ and $S=0$.
\item[Case I.c] $s=x_1^3$ and $S=X_3^3$: $s$ is a triple line $\ell^3$ in
$\P V^*$ and $S$ is a triple line $p^3$ in $\P V$, with $p$ (viewed as a
point in $\P V^*$) lying on $\ell$.
\item[Case I.d] $s=x_1^2x_2$ and $S=X_2X_3^2$: $s=\ell^2\cup\ell'$ and
$S=p^2\cup p'$, with $p=\ell\cap\ell'$ and $p'\in\ell$.
\item[Case I.e] $s=\al\,x_1x_2x_3$ and $S=\al\,X_1X_2X_3$, with
$\al\ne0,1,-1$: $s=\ell\cup\ell'\cup\ell''$ and $S=p\cup p'\cup p''$, with
$\ell\cap\ell'\cap\ell''=\emptyset$, $p=\ell'\cap\ell'$,
$p'=\ell''\cap\ell$, $p''=\ell\cap\ell'$.
\item[Case I.f] $s=6i\sqrt3\,x_1(x_1^2+x_2x_3)$ and
$S=6i\sqrt3\,X_1X_2X_3$, where $i^2=-1$: same configuration as case I.e,
but $s=\Con\cup\ell$ and $S=p\cup p'\cup p''$, with $\Con$ a
(nondegenerate) conic tangent to $\ell'$ at $p''$ and to $\ell''$ at $p'$.
\item[Case I.g] $s=6i\sqrt3\,x_1(x_1^2+x_2x_3)$ and
$S=6i\sqrt3\,X_3(X_3^2+X_1X_2)$: same configuration as case I.f, but
$s=\Con\cup\ell$ and $S=\Con'\cup p'$, with $\Con'$ a conic tangent to $p$
at $\ell''$ and to $p''$ at $\ell$.
\item[Case I.h] $e=\al(x_{123}+x_{231}+x_{312})
+\be(x_{132}+x_{213}+x_{321}) +\ga(x_{111}+x_{222}+x_{333})$ and
$E=\al'(x^{123}+x^{231}+x^{312}) +\be'(x^{132}+x^{213}+x^{321})
+\ga'(x^{111}+x^{222}+x^{333})$, with $\ga\ne0$, $\ga'\ne0$,
$\ga^3+(\al+\be)^3\ne0$, ${\ga'}^3+(\al'+\be')^3\ne0$ (so $s$ and $S$ are
elliptic curves), $\al\ne\be$, $\al'\ne\be'$, $\al\al'+\be\be'+\ga\ga'\ne0$
and
\begin{equation}
\label{ellcond}
\left\{
\begin{gathered}
\al^2{\al'}^2+\be^2{\be'}^2+\ga^2{\ga'}^2
-2\al\al'\be\be'-2\al\al'\ga\ga'-2\be\be'\ga\ga'=0\\
\al^2\be'\ga'+\be^2\al'\ga'+\ga^2\al'\be'=0\\
{\al'}^2\be\ga+{\be'}^2\al\ga+{\ga'}^2\al\be=0
\end{gathered}
\right.
\end{equation}
(In this case, we refrain from normalizing the antisymmetric parts
$\al-\be$ and $\al'-\be'$, to keep \eqref{ellcond} homogeneous.) Note that
there are solutions, e.g., $\al=\al'=0$ and $\be=\be'=\ga=\ga'\ne0$.

{\bf Question.} Can conditions \eqref{ellcond} be described geometrically
in terms of the elliptic curves $s$ and $S$?
\end{description}
(The cases where $s$ is a cusp curve, a node curve, a conic with a tangent
line or three intersecting lines cannot occur.)

The case $\la\ne0$, $\La=0$ (or vice-versa) is impossible.

There is only one possible case with $\la=0$, $\La=0$.
\begin{description}
\item[Case I.e$^*$] $\la=0$, $\La=0$, $s=x_1x_2x_3$ and $S=X_1X_2X_3$.
(The geometry is similar to case I.e.)
\end{description}

\subsection{Type~I'}

We take $C,D,c,d$ as for type I. Define
\[
z_{ijk}:=x_{ijk}+\om x_{jki}+\om^2 x_{kij}
\qquad
z^{ijk}:=x^{ijk}+\om x^{kij}+\om^2 x^{jki}
\]
Let $\Gamma_\om$ be the subspace of $V^{\o3}$ spanned by all $z_{ijk}$ and
$\Gamma_\om^*$ that of ${V^*}^{\o3}$ spanned by all $z^{ijk}$. Now
\eqref{QE} means that $e\in\Gamma_\om$ and $E\in\Gamma^*_\om$.

Note that $\Gamma_\om$ is a sub-$\GL(V)$-module of $V^{\o3}$, isomorphic to
the $\GL(V)$-module $\sll(V)=\{X\in\Lin(V,V)\mid\tr X=0\}$. We use the
isomorphism given by
\[
\pmat{t_1&u_1&u_3\\v_1&t_2-t_1&u_2\\v_3&v_1&-t_2}\mapsto
\begin{aligned}
{}&t_1(z_{123}+z_{213})+t_2(z_{231}+z_{321})\\
{}&\quad-u_1z_{113}-u_2z_{221}+u_3z_{112}\\
{}&\quad+v_1z_{223}+v_2z_{331}-v_3z_{332}
\end{aligned}
\]
to view $e$ as an element of $\sll(V)$. We do the same for $E$, using the
isomorphism
\[
\pmat{t'_1&u'_1&u'_3\\v'_1&t'_2-t'_1&u'_2\\v'_3&v'_1&-t'_2}\mapsto
\begin{aligned}
{}&t'_1(z^{123}+z^{213})+t'_2(z^{231}+z^{321})\\
{}&\quad-u'_1z^{223}-u'_2z^{331}+u'_3z^{332}\\
{}&\quad+v'_1z^{113}+v'_2z^{221}-v'_3z^{112}
\end{aligned}
\]
between $\Gamma^*_\om$ and $\sll(V)$.

Using $\Stab(Q)=\GL(V)$, we reduce $e$ to a Jordan normal form. Since $\tr
e=0$ and $e\ne0$, there are two cases.
\begin{description}
\item[Case I'.a] $e$ is diagonal, i.e. (after rescaling), $u_i=v_i=0$
($i=1,2,3$), $t_1=1$, $t_2\ne2$, $t_2\ne\frac12$. Now \eqref{Ee} reads
\begin{gather*}
u'_1=u'_2=v'_1=v'_2=0\qquad u'_3(t_2+1)=v'_3(t_2+1)=0\\
(t_2-2)t'_1+(1-2t_2)t'_2\ne0
\end{gather*}
There are two subcases.
\begin{itemize}
\item If $t_2\ne-1$, then $u'_3=v'_3=0$, and (\ref{coh}gh) are equivalent
to $t_2t'_1=t'_2$. Rescaling $e$, we get $t'_1=1$, $t'_2=t_2$.
\item If $t_2=-1$, then (\ref{coh}gh) are equivalent to $u'_3=v'_3=0$,
$t'_2=-t'_1$. Rescaling $E$, we get $t'_1=1$, $t'_2=-1$.
\end{itemize}
Combining these subcases gives a $1$-parameter family
\[
e=\pmat{1&0&0\\0&t-1&0\\0&0&-t}
\qquad
E=\pmat{1&0&0\\0&t-1&0\\0&0&-t}
\]
with the condition $t^2-t+1\ne0$.
\item[Case I'.b] $e$ has a $2\times2$ Jordan block, i.e. (after
rescaling), $v_1=v_2=v_3=0$, $u_1=u_2=0$, $u_3=1$, $t_1=-t_2=1$. Now
\eqref{Ee} reads
\[
v'_1=v'_2=v'_3=0\qquad u'_1=u'_2=0\qquad t'_1=-t'_2\ne0
\]
and (\ref{coh}gh) are equivalent to $t'_1=-u'_3$. Rescaling $E$, we obtain
the solution
\[
e=\pmat{1&0&1\\0&-2&0\\0&0&1}
\qquad
E=\pmat{1&0&-1\\0&-2&0\\0&0&1}
\]
\end{description}

\subsection{Type~II}

We take
\[
\begin{aligned}
c&=q\inv\,x_1\o y_1+x_2\o y_2+q\,x_3\o y_3\\
C&=q\,y^1\o x^1+y^2\o x^2+q\inv\,y^3\o x^3
\end{aligned}
\qquad
\begin{aligned}
d&=q\,y_1\o x_1+y_2\o x_2+q\inv\,y_3\o x_3\\
D&=q\inv\,x^1\o y^1+x^2\o y^2+q\,x^3\o y^3
\end{aligned}
\]
By \eqref{QE}, $e$ and $E$ are of the form
\begin{align*}
e&=\al(x_{123}+q^2\,x_{231}+q^2\,x_{312})
-\be(x_{132}+x_{213}+q^2\,x_{321})
+\ga x_{222}\\
E&=\al'(x^{123}+q^2\,x^{231}+q^2\,x^{312})
-\be'(x^{132}+x^{213}+q^2\,x^{321})
+\ga'x^{222}
\end{align*}
Condition \eqref{Ee} then reads
\[
(q^2-1)(q^2\al\al'-\be\be')+\ga\ga'=0\qquad\al\al'+\be\be'\ne0
\]
Inspecting $\text{(\ref{coh}g)}-\tp\text{(\ref{coh}h)}$ shows that we must
have $\ga\ga'=0$, so up to Dynkin flip (cf.\ section~\ref{equivsection}),
we may assume $\ga'=0$. Now (\ref{coh}gh) become equivalent to
$\ga(q^4\al^3-\be^3)=0$. We therefore have two cases.
\begin{description}
\item[Case II.a] If $\ga=0$, we normalize $e,E$ so that $\al=\al'=1$.
This gives a $2$-parameter family, namely the Artin-Schelter-Tate quantum
$SL(3)$'s~\cite{AST} (or rather, their quantum $\GL(3)$'s having a central
quantum determinant), where (in the notation of~\cite{AST})
$p_{21}=p_{32}=\be$, $p_{31}=\be'$ and $\la=\frac1{\be\be'}=\frac1{q^2}$.
The standard quantum $\SL(3)$~\cite{FRT} is obtained as a particular case,
when $\be=\be'=q$.
\item[Case II.b] If $\ga\ne0$, then $\be=p^4\al$, with $p^3=q$. We
first normalize $e,E$ so that $\al=p^{-2}$, $\be=p^2$, $\al'=p^{-1}$ and
$\be'=p$, then we use $\Stab(Q)=\{\spmat{\la&0&0\\0&\mu&0\\0&0&\nu}\}$ to
get $\ga=q^2-1$. This gives a $1$-parameter family, namely the
Cremmer-Gervais quantum $\SL(3)$ as described in~\cite{Ho}.
\end{description}

\subsection{Type~II'}

We take $C,D,c,d$ as for type II. By \eqref{QE}, $e$ and $E$ are of the
form
\begin{align*}
e&=\al(x_{123}+\om q^2\,x_{231}+\om^2q^2\,x_{312})
-\be(x_{132}+\om q^2\,x_{321}+\om^2\,x_{213})\\
E&=\al'(x^{123}+\om q^2\,x^{312}+\om^2q^2\,x^{231})
-\be'(x^{132}+\om\,x^{213}+\om^2q^2\,x^{321})
\end{align*}
Condition \eqref{Ee} then reads
\[
\be\be'=q^2\al\al'\qquad \al\al'+\be\be'\ne0
\]
\begin{description}
\item[Case II'.a] (unique case of this type) Rescaling $e,E$, we get
$\al=\al'=1$ and $\be'=\frac{q^2}{\be}$. Since conditions (\ref{coh}gh) are
automatically fulfilled, this gives a $2$-parameter family.
\end{description}

\subsection{Type~III}

We take
\begin{align*}
c&=x_1\o y_1+x_2\o y_2+x_3\o y_3+x_1\o y_3\\
C&=y^1\o x^1+y^2\o x^2+y^3\o x^3-y^1\o x^3\\
d&=y_1\o x_1+y_2\o x_2+y_3\o x_3-y_3\o x_1\\
D&=x^1\o y^1+x^2\o y^2+x^3\o y^3+x^3\o y^1
\end{align*}
By \eqref{QE}, $e$ and $E$ are of the form
\begin{align*}
e&=\al(x_{123}+x_{231}+x_{312}-x_{321}-x_{213}-x_{132})-2\al\,x_{121}\\
&\quad+\be\,x_{111}+\ga(x_{112}+x_{121}+x_{211})
+\de(x_{122}+x_{212}+x_{221})+\eps\,x_{222}\\
E&=\al'(x^{321}+x^{213}+x^{132}-x^{123}-x^{231}-x^{312})-2\al'x^{323}\\
&\quad+\be'x^{333}+\ga'(x^{332}+x^{323}+x^{233})
+\de'(x^{322}+x^{232}+x^{223})+\eps'x^{222}
\end{align*}
Now \eqref{Ee} reads
\[
\de\de'=\de\eps'=\eps\de'=\eps\eps'=0\qquad\al\al'\ne0
\]
We normalize to $\al=\al'=1$.

If $\eps=\eps'=0$, then (up to Dynkin flip) $\de'=0$. Now (\ref{coh}gh) is
equivalent to
\[
\de=0\qquad(\ga+\ga'-1)(\ga+\ga'-2)=0
\]
Note also that if $\ga\ne\frac23$ (resp.\ $\ga'\ne\frac23$), then we may
use $\Stab(Q)=\{\spmat{\la&\sig&\nu\\0&\mu&\sig'\\0&0&\la}\}$ to get
$\be=0$ (resp.\ $\be'=0$). Working up to $\Stab(Q)$ (and up to Dynkin flip)
now leads to four cases (two of which are $1$-parameter families).
\begin{description}
\item[Case III.a] $(\al,\be,\ga,\de,\eps)=(1,0,\ga,0,0)$ and
$(\al',\be',\ga',\de',\eps')=(1,0,1-\ga,0,0)$.
\item[Case III.a$^*$] $(\al,\be,\ga,\de,\eps)=(1,1,\frac23,0,0)$ and
$(\al',\be',\ga',\de',\eps')=(1,0,\frac13,0,0)$.
\item[Case III.b] $(\al,\be,\ga,\de,\eps)=(1,0,\ga,0,0)$ and
$(\al',\be',\ga',\de',\eps')=(1,0,2-\ga,0,0)$.
\item[Case III.b$^*$] $(\al,\be,\ga,\de,\eps)=(1,1,\frac23,0,0)$ and
$(\al',\be',\ga',\de',\eps')=(1,0,\frac43,0,0)$.
\end{description}

If one of $\eps,\eps'$ is nonzero, then up to Dynkin flip, we may assume
$\eps\ne0$. If follows that $\de'=\eps'=0$. Using $\Stab(Q)$, we may get
$\eps=1$ and $\de=0$. Now (\ref{coh}gh) is equivalent to
\[
\be=0\qquad \ga'=\frac12(1+\ga)\qquad(\ga+\ga'-1)(\ga+\ga'-2)=0
\]
This leads to two further cases.
\begin{description}
\item[Case III.c] If $\ga=\frac13$, we have a $1$-parameter family
\[
(\al,\be,\ga,\de,\eps)=(1,0,\frac13,0,1)
\qquad
(\al',\be',\ga',\de',\eps)=(1,\be',\frac23,0,0)
\]
\item[Case III.c$^*$] If $\ga=1$, we may still use $\Stab(Q)$ to get
$\be'=0$.
This gives the solution
\[
(\al,\be,\ga,\de,\eps)=(1,0,1,0,1)
\qquad
(\al',\be',\ga',\de',\eps')=(1,0,1,0,0)
\]
\end{description}

\subsection{Type~III'}

We take $C,D,c,d$ as for type III. By \eqref{QE}, $e$ and $E$ are of the
form
\begin{align*}
e&=\al(x_{123}+\om\,x_{231}+\om^2\,x_{312}
-x_{321}-\om\,x_{213}-\om^2\,x_{132})
-2\al\,x_{121}\\
&\quad+\be\,x_{111}
+\frac{\om-1}2\be(x_{113}+\om\,x_{131}+\om^2\,x_{311})\\
&\quad+\ga(x_{112}+\om\,x_{121}+\om^2\,x_{211})
+\de(x_{122}+\om\,x_{221}+\om^2\,x_{212})
&\quad\\
E&=\al'(x^{321}+\om^2\,x^{213}+\om\,x^{132}
-x^{123}-\om^2\,x^{231}-\om\,x^{312})
-2\al'x^{323}\\
&\quad+\be'x^{333}
+\frac{\om^2-1}2\be'(x^{331}+\om^2\,x^{313}+\om\,x^{133})\\
&\quad+\ga'(x^{332}+\om^2\,x^{323}+\om\,x^{233})
+\de'(x^{322}+\om^2\,x^{223}+\om\,x^{232})
\end{align*}
Now \eqref{Ee} reads
\[
\begin{aligned}
3\al\be'&=2(1-\om)\de\al'\\
3\al'\be&=2(1-\om^2)\de'\al
\end{aligned}
\qquad
\begin{aligned}
3\be\be'&=4\de\de'\\
\al\al'&\ne0
\end{aligned}
\]
It follows that $\de\de'=0$, so up to Dynkin flip, we may assume that
$\de'=0$. Rescaling $e,E$ so that $\al=\al'=1$, we get $\be=0$. Using
$\Stab(Q)$, we may get $\de=0$, so $\be'=0$. Now (\ref{coh}gh) are
equivalent to
\[
(\ga+\om\ga'-\om^2)(\ga+\om\ga'-2\om^2)=0
\]
We therefore have two $1$-parameter families.
\begin{description}
\item[Case III'.a] $(\al,\be,\ga,\de)=(1,0,\ga,0)$ and
$(\al',\be',\ga',\de')=(1,0,\om-\om^2\ga,0)$
\item[Case III'.b] $(\al,\be,\ga,\de)=(1,0,\ga,0)$ and
$(\al',\be',\ga',\de')=(1,0,2\om-\om^2\ga,0)$
\end{description}

\subsection{Type~IV}

We take
\begin{align*}
c&=x_1\o y_1+x_2\o y_2+x_3\o y_3+x_1\o y_2+x_2\o y_3+\tfrac12\,x_1\o y_3\\
C&=y^1\o x^1+y^2\o x^2+y^3\o x^3-y^1\o x^2-y^2\o x^3+\tfrac12\,y^1\o x^3\\
d&=y_1\o x_1+y_2\o x_2+y_3\o x_3-y_2\o x_1-y_3\o x_2+\tfrac12\,y_3\o x_1\\
D&=x^1\o y^1+x^2\o y^2+x^3\o y^3+x^2\o y^1+x^3\o y^2+\tfrac12\,x^3\o y^1
\end{align*}
By \eqref{QE}, $e$ and $E$ are of the form
\begin{align*}
e&=\al(x_{123}+x_{231}+x_{312}-x_{321}-x_{213}-x_{132})\\
&\quad+2\al(x_{112}-x_{211})-2\al(x_{113}+x_{311})+2\al\,x_{212}\\
&\quad+\be\,x_{111}+\de(x_{122}+x_{212}+x_{221})\\
&\quad+2\de(x_{112}-x_{211})-2\de(x_{113}+x_{131}+x_{311})\\
E&=\al'(x^{321}+x^{213}+x^{132}-x^{123}-x^{231}-x^{312})\\
&\quad+2\al'(x^{332}-x^{233})-2\al'(x^{331}+x^{133})+2\al'x^{232}\\
&\quad+\be'x^{333}+\de'(x^{322}+x^{232}+x^{223})\\
&\quad+2\de'(x^{332}-x^{233})-2\de'(x^{331}+x^{313}+x^{133})
\end{align*}
Now \eqref{Ee} reads
\[
2(\al\de'+\al'\de)+9\de\de'=0\qquad\al\al'\ne0
\]
We normalize to $\al=\al'=1$. Next,
$\text{(\ref{coh}g)}-\tp\text{(\ref{coh}h)}$ reads
$(\de'-2\de)(\de-2\de')=0$, so up to Dynkin flip, we may assume that
$\de'=2\de$. We now have two cases.
\begin{description}
\item[Case IV.a] If $\de=0$, then (\ref{coh}gh) are equivalent to
$\be+\be'+2=0$. Using
$\Stab(Q)=\{\spmat{\la&\mu&\nu\\0&\la&\mu\\0&0&\la}\}$, we may get
$\be=\be'=-1$; this gives the solution
\[
(\al,\be,\de)=(1,-1,0)\qquad(\al',\be',\de')=(1,-1,0)
\]
\item[Case IV.b] If $\de=-\frac13$, then (\ref{coh}gh) are equivalent
to
$\be'=-\frac8{27}$. Using $\Stab(Q)$, we may get $\be=0$; this gives the
solution
\[
(\al,\be,\de)=(1,0,-\frac13)\qquad(\al',\be',\de')=(1,-\frac8{27},-\frac23)
\]
\end{description}

\subsection{Type~IV'}

We take $C,D,c,d$ as for type IV. Condition \eqref{QE} prevents $e,E$ from
satisfying \eqref{Ee}, so this type is impossible.

\section{Further problems}
\label{outlooksection}

For $G=\SL(3)$, some technical problems and some links with other
literature should be worthwhile studying:
\begin{itemize}
\item Replace the case by case argument in the proof of
proposition~\ref{Mright} by a more conceptual one, preferrably including
the elliptic case.
\item Classify the elliptic solutions more explicitely, i.e., study
conditions \eqref{ellcond} (which define a subvariety in $\P^2\times\P^2$).
\item View each quantum $\SL(3)$ with $\om=1$ as a formal $1$-parameter
deformation of $\SL(3)$, compute the Lie bialgebra structure on $\sll(3)$
at its semi-classical limit (see, e.g., \cite{CP} for definitions) and
compare with the classification of these structures given in~\cite{St}.
(Obviously, a quantum $\SL(3)$ with $\om\ne1$ cannot be viewed as such a
deformation.) Note that a related converse problem---that of finding an
$R$-matrix quantizing each of the Lie bialgebra structures on
$\sll(3)$---has recently been solved in~\cite{GG}.
\item Determine which quantum $\SL(3)$'s admit a compact form. (The
standard quantum $\SL(2)$ has a compact form, but the Jordanian one does
not~\cite{Ohn}.)
\item Our definition of a BQD seems to be related to the spiders
of type A$_2$ introduced in~\cite{Ku}.
\item The shape algebras $\M_\L$ and $\MM_\L$ from
section~\ref{shapesection} being homogeneous quadratic algebras, try to use
the methods of~\cite{ATV} to associate geometric data to such a pair of
algebras, in order to get a better understanding of the classification of
quantum $\SL(3)$'s. Incidentally, we note that the matrix $Q$ used in
section~\ref{classsection} is the same as that used in~\cite{AS} to
classify regular algebras of dimension~$3$ (here, such algebras would arise
as quantum analogues of the homogeneous coordinate ring of $\SL(3)/P$, $P$
a maximal parabolic subgroup).
\end{itemize}
Another problem is of course to extend the methods used here to study
quantum $G$'s for some other reductive group $G$. However, a classification
is probably out of reach, even for $G=\SL(4)$.

\end{document}